\documentclass[a4paper,11pt]{article}
\usepackage{latexsym}
\usepackage{graphicx}
\usepackage{hyperref}
\usepackage{color}
\usepackage{amssymb,amsmath}
\usepackage{multicol}
\usepackage{cancel}
\usepackage{multirow}
\usepackage[utf8]{inputenc}
\usepackage[T1]{fontenc}
\usepackage{authblk}
\usepackage{epsfig}
\usepackage{epstopdf}
\usepackage{caption}
\usepackage{multirow}
\usepackage{adjustbox}
\usepackage{array}
\usepackage{tabularray}

 \DeclareMathOperator{\erf}{erf}
 \newtheorem{exm}{Example}
  \newtheorem{alg}{Algorithm}
\title{Target Detection in a Known Number of Intervals Based on Cooperative Search Technique}
\author[1,2]{M. Fakharany\thanks{fakharany@aucegypt.edu (mohamed.elfakharany@science.tanta.edu.eg)}}
\author[3]{Faten Alamri \thanks{fsalamri@pnu.edu.sa}}
\author[1,2]{Mohamed Abd Allah El-Hadidy \thanks{melhadidi@science.tanta.edu.eg}}
\affil[1]{Mathematics Department, Faculty of Science, Tanta University, Tanta, Egypt.}
\affil[2]{Mathematics and Statistics Department, College of Science, Taibah University, Yanbu, Saudi Arabia.}
\affil[3]{Mathematical Sciences Department, College of Science, Princess Nourah Bint Abdulrahman University,
P.O. Box 84428, Riyadh 11671, Saudi Arabia.}

\begin{document}
  \maketitle
  \begin{abstract}
Finding hidden/lost targets in a broad region costs strenuous effort and takes a long time. From a practical view, it is convenient to analyze the available data to exclude some parts of the search region. This paper discusses the coordinated search technique of a one-dimensional problem with a search region consisting of several mutual intervals. In other words, if the lost target has a probability of existing in a bounded interval, then the successive bounded interval has a far-fetched probability. Moreover, the search domain is swept by two searchers moving in opposite directions, leading to three categories of target distribution truncations; commensurate, uneven, and symmetric. The truncated probability distributions are defined and applied based on the proposed classification to calculate the expected value of the elapsed time to find the hidden object. Furthermore, the optimization of the associated expected time values of various cases is investigated based on Newton's method. Several examples are presented to discuss the behavior of various distributions under each case of truncation. Also, the associated expected time values are calculated as their minimum values.     
\end{abstract}
\section{Introduction}
Usually, searching for lost objects in the sea/ocean or land requires a lot of effort, money, and advanced devices. The lost object in the ocean can be a person, ship, airplane, or submarine, while on the land, it can be antiquities or mine. Moreover, when an object is lost in the sea, it suffers from the drift of water and air currents. For instance, the Malaysian airlines Flight 370 was found close to the eastern coast of Africa in the Indian Ocean. This region is so far from the best guess location by thousands of kilometers.\\\\
The search theory provides mathematical models and approaches to optimize the search process. \cite{koopman1946search} and \cite{ stone1977search} are two pioneers who contributed to developing this field, especially in anti-submarine warfare during World War II. Also, this problem is considered in one, two, and three dimensions, where the search domain is line, plane, and space, respectively. Since then, various mathematical models have enhanced the search process, improving the best guess location and optimizing the time required.\\\\
From a practical viewpoint, the search domain is required to truncate. Usually, truncation is applied to convert the unbounded search zone into a bounded one. Typically, planners exploit all the information about the target and then determine the probability of its best guess location. In this case, the probability of other places is not considered. Thus, we need to truncate the probability distribution of these neglected places and then distribute it to the searching places (places with the highest probability). Consequently, it helps us to detect the target as soon as feasible.\\\\
The coordinated search technique is one of the famous techniques used to detect the target as quickly as possible. \cite{thomas1992finding} provided the first investigation of this technique to detect a hidden target on the circle circumference. He calculated the expected value of the detection time after the target was detected. On the real line, \cite{reyniers1996coordinated} discussed this technique using two searchers with equal speeds. These searchers start the searching process from the origin of the line. Each searcher specializes in its part (left or right) search. After finishing each step of the searching process, the searchers will return to the origin. Also, \cite{reyniers1995co} applied these tools to study this technique on the known interval to detect the hidden target. More than obtaining the expected detection time value, she derived the optimal search technique to minimize this time. \cite{alfreedi2019optimal} and  \cite{el2022linear} studied this technique on the line when the target is randomly located with known distribution (symmetric or asymmetric). They obtained the necessary conditions for optimal search technique and, consequently, the minimum expected value of the detection time obtained. In the case of the randomly moving target (i.e., when the target moves with a known stochastic process such as Random Walk and Brownian motion processes), \cite{el2018coordinated, el2020existence, el2019cooperative,el2020detection,el2016searching,el2021optimal11,Abdallah2022study12,el2021study13,el2021existence14} discussed different analytical studies to show the existence of the finite search technique. More than obtaining the necessary conditions for optimality of the search technique, they presented the minimum expected value of the first meeting time between one of the searchers and the target. In addition, this technique was discussed on the plane and space using multiple searchers as in  \cite{abou2009coordinated15,el2012coordinated16,el20183ref17,teamah2021probabilistic18}. They studied the optimality and the existence of the finiteness of this technique. For more studies on the search theory for missing targets, \cite{el2016fuzzy22,mohamed2015optimal23,mohamed2011multiplicative24,el2019generalised25,el2020existence26,el2019existence27,el2019mathematical28,
el2019searching29,el2019studying30,alzulaibani2019study31,el2021existence32,kassem2014optimal33,el2021stochastically34,el2021developing35,
el2021optimal36,el2021optimal37,el2022detection38,el2021quality39} provided more statistical and analytical studies for different search techniques such as the linear search technique.\\\\
So far, the truncation is applied on one of the endpoints of the search zone or both. However, it can be substantially vast, especially thousands of kilometers. Consequently, the search becomes impractical. This paper investigates several truncations within a truncated search zone. We consider a search problem of a hidden target following some distributions in a one-dimensional problem. The proposed model consists of two searchers sweeping a line of finite length with deleted subintervals. They start the searching process from the origin (the observer center) of the real line; the first sweeps towards the right and the other towards the left. Several subintervals are truncated in each part of the line based on the collected information. Consequently, they skip these zones. Applying several truncations leads to the truncation type based on the numbers in each part of the line and whether they are symmetrical. Consequently, we establish a classification of the truncations. Next, a truncated distribution is defined on a domain with deleted subintervals. The expected elapsed time to find the lost target is determined for each type of truncation, which is a function of the mesh truncation points. An optimization technique has been applied to obtain the best value of expectation. Also, we consider a particular case when one searcher scans a truncated half-line with some deleted subintervals.   \\\\
The organization of this paper appears as follows. Section 2 introduces the truncated subintervals within a search zone, and the classification of the used truncations is established. The probability distribution functions and their cumulative distribution functions are defined for each case of truncation. Also, the mechanism of the search process is introduced. The expected elapsed time value is introduced in Section 3 for each category of the applied truncation. Next, the optimization of its value is investigated. Section 4 is dedicated to implementing several computational examples to reveal the properties of truncated distributions and calculate the associated expected elapsed time values for several truncated distributions, and obtain the corresponding optimized values. The conclusion and future work are discussed in the last section.\\
\section{The mechanism of the searchers}
\begin{figure}[h]
  \hspace{-0.0cm}\includegraphics[width=1.3\linewidth]{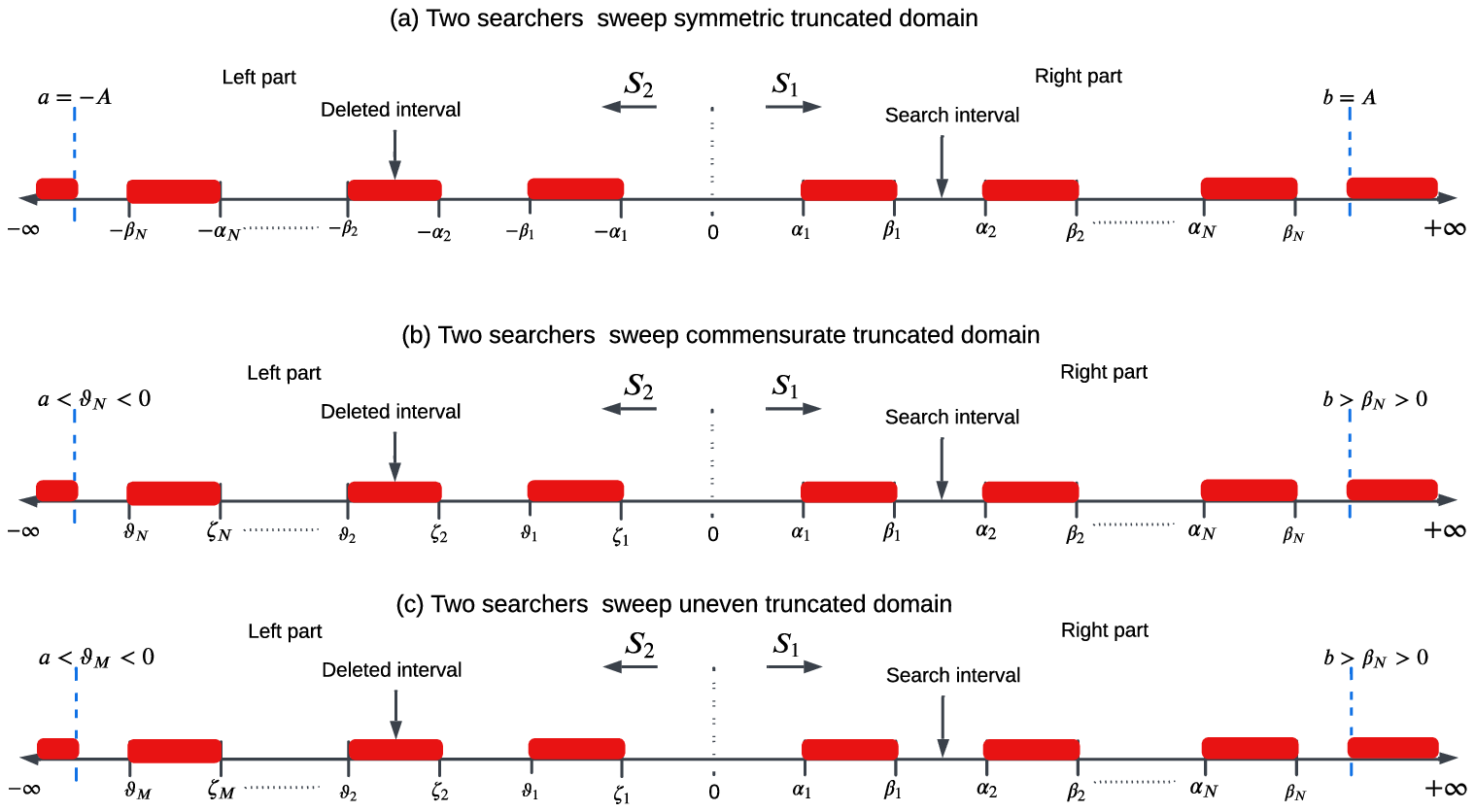}
  \vspace{-1.5cm}\caption{\footnotesize The search space (real line) based on various type of truncations.}\label{Fig1}
\end{figure}

Naturally, there is some available information about the target's position. This information can determine the presence of the target in some places. Instead of searching for the target in a whole line, as happened in the previous works, see, for example, \cite{reyniers1995co} and \cite{alfreedi2019optimal} and  \cite{el2022linear}, we use this information to determine the most likely places for the target to be located. Thus, we can consider the searching space as a set of known intervals on the real line. These intervals will be determined according to the available information about the target's position and the probability of the target being present in each interval. The probability of the target in the intervals which have no information (deleted intervals from the line) is incorporated into the probability of the searching intervals. To do this, we use the mathematical definition presented in \cite{el2019generalized40} to delete a set of intervals from the distribution domain. Suppose that $X$ is a random variable with a probability density function (PDF) $f_{X}(x)$ of the target location with cumulative distribution function (CDF) $F_{X}(x)$, and $x\in\mathbb{R}=(-\infty,+\infty)$. Based on the provided data, we guarantee that the lost target lies in an interval, say $I=[a,b]$, moreover, it cannot be in some subintervals within the interval $I$. Consider a set $\mathcal{P}$ of  finite sequence of mesh-points, given by
\begin{equation}\label{PertPointsOne}
\begin{array}{rl}
  \mathcal{P}= & \{a<\vartheta_{M}<\zeta_{M}<\vartheta_{M-1}<\zeta_{M-1}<\ldots<
  \vartheta_{1}<\zeta_{1}<0 \\
   & <\alpha_{1}<\beta_{1}<\alpha_{2}<\beta_{2}<\ldots<\alpha_{N}<\beta_{N}<b\}.
\end{array}
\end{equation}
The open deleted subintervals; $\{\breve{I}_{j}=(\vartheta_{j},\zeta_{j})\}_{j=1}^{M}$ for the left part of the real line and $\{\hat{I}_{j}=(\alpha_{j},\beta_{j})\}_{j=1}^{M}$ for its right part, represent the intervals where the lost target cannot exist in them.  In other words, the unbounded domain $\mathbb{R}$ is contracted into the interval $I=[a,b]$, and we have the deleted subintervals $\{\breve{I}_{j}\}_{j=1}^{N}$, and $\{\hat{I}_{j}\}_{j=1}^{N}$ for the negative and positive values of the truncated real line. Consequently, we can classify the truncations into three categories based on the value of $M$, and  $\{\breve{I}_{j}\}_{j=1}^{M}$. First category when $M=N$, $b=-a=A\in\mathbb{R}^{+}$, and $\beta_{j}=-\vartheta_{j}$, and $\alpha_{j}=-\zeta_{j}$, for all values of $j\in\{1,2,\dots,N\}$, then we have symmetric truncation. Second, $M=N$, but $a\neq b$, $\beta_{j}\neq-\vartheta_{j}$, or $\alpha_{j}\neq-\zeta_{j}$ for at least one value of $j$, then this truncation is called commensurate truncation. Finally, for $M\neq N$, the truncation is said to be uneven truncation. In light of this classification, we use $\tilde{M}$ to denote the number of deleted subinterval in the negative part of the real line such that $\tilde{M}=N$, for the symmetric and commensurate truncations, while $\tilde{M}=M$, for the uneven truncation. Let $\hat{N}$ be the total number of the truncated subintervals, the $\hat{N}=2N$ in the commensurate truncation, while $\hat{N}=N+M$, for the uneven case. Therefore, we have the compact domain $\Omega$, given by
\begin{equation}\label{Domaineq1}
\begin{array}{l}
    \Omega=\left(\displaystyle{\bigcup_{j=1}^{\tilde{M}}}[\zeta_{j+1},\vartheta_{j}]\right)
\cup[\zeta_{1},\alpha_{1}]\cup\left(\displaystyle{\bigcup_{\ell=1}^{N}}[\beta_{\ell},\alpha_{\ell+1}]\right)=\Omega^{-}\cup\Omega^{+}, \\
    \Omega^{-}=\left(\displaystyle{\bigcup_{j=1}^{\tilde{M}}}[\zeta_{j+1},\vartheta_{j}]\right)
\cup[\zeta_{1},0], \\
  \Omega^{+}=[0,\alpha_{1}]\cup\left(\displaystyle{\bigcup_{\ell=1}^{N}}[\beta_{\ell},\alpha_{\ell+1}]\right),
  \end{array}
\end{equation}
where $\zeta_{\tilde{M}+1}=a$, and $\alpha_{N+1}=b$. Here, we have two searchers $S_{1}$ and $S_{2}$ sweep $\Omega^{+}$ and $\Omega^{-}$, respectively, and they are connected with the observer center to update the sweeping results at once. Each of them starts from the origin, $S_{1}$ moves towards the right with speed $v_{1}$ and $S_{2}$ moves with the same speed towards the left as shown in Figure \ref{Fig1}. Moreover, both of them move with speed $v_{2}$ in the deleted subintervals ($v_{2}>v_{1}$).\\
Before, we investigate the truncated distributions, we introduce a synoptic view of distributions and their properties. First, the distribution is said to symmetric if and only if there exists a number $\hat{x}$ such that
$$f_{X}(\hat{x}-\xi)=f_{X}(\hat{x}+\xi),~~\forall\xi\in\mathbb{R},$$
otherwise the distribution is called skew-symmetric (asymmetric). Next the domain of the distribution can be bounded such as in beta distribution, bounded above, bounded below, or unbounded domain. In this paper, we interest in symmetric and asymmetric distributions defined on unbounded domain $(-\infty,+\infty)$. Also, we consider the distributions defined on bounded below domain $[0,+\infty)$. The truncations of these domains are introduced, analyzed, and applied on some distributions that listed in Table \ref{Tab1}, also, their CDFs are given in Table \ref{Tab2}. The parameters $\mu$ and $\sigma$ are the mean, and standard deviation of the normal distribution, $\tilde{x}\in\mathbb{R}$ represents the location of the Cauchy distribution peak and $c>0$ is half-width at half-maximum of the distribution, $\kappa,$ $\theta>0$ are the shape-scale parameters of gamma distribution. The skew normal distribution has the location, shape, and scale parameters $\eta$, $\varrho\in\mathbb{R}$, and $\varpi>0$. The functions $\Phi(x)$, and $\mathcal{T}(x,\tilde{a})$ are defined by
\begin{equation}\label{SkewNM1}
  \begin{array}{rl}
    \Phi(x)=&\frac{1}{2}\left(1+\erf\left(\frac{x}{\sqrt{2}}\right)\right),\\
    \mathcal{T}(x,\tilde{a})=& \frac{1}{2\pi}\displaystyle{\int_{0}^{\tilde{a}}}\frac{e^{-\frac{x^{2}}{2}(1+\nu^{2})}}{1+\nu^{2}}d\nu,~~x,\tilde{a}\in\mathbb{R},
  \end{array}
\end{equation}
where, $\mathcal{T}(x,\tilde{a})$ is the Owen's T-function \cite{owen1956tables}. The special functions $\erf(x)$, $\Gamma(x)$, and $\gamma(\kappa,x)$ are the error, gamma, and lower incomplete gamma functions.
\begin{center}
\begin{table}
\begin{center}
  \begin{tabular}{p{4.5cm}p{6.5cm}p{3.0cm}}
    \hline
    \textbf{Distribution name} & \textbf{The PDF} & \textbf{Symmetric/}\\
     & & \textbf{Asymmetric}\\
    \hline
    Normal distribution & $\displaystyle{f(x)=\frac{1}{\sigma\sqrt{2\pi}}e^{\frac{-1}{2}\left(\frac{x-\mu}{\sigma}\right)^{2}}}$ & Symmetric\\
    Cauchy distribution & $\displaystyle{f(x)=\frac{c}{\pi(c^{2}+(x-\tilde{x})^{2})}}$& Symmetric  \\
    Skew-normal distribution &$\displaystyle{f(x)=\frac{1}{\varpi\sqrt{2\pi}}e^{-\frac{(x-\eta)^{2}}{2\varpi^{2}}}
    \Phi\left(\frac{\varrho}{\varpi}(x-\eta)\right)}$  & Asymmetric \\
    Gamma distribution & $\displaystyle{f(x)=\frac{x^{\kappa-1}e^{-\frac{x}{\theta}}}{\Gamma(\kappa)\theta^{\kappa}}}$& Asymmetric \\
    \hline
  \end{tabular}
\end{center}
 \caption{The PDF and properties of some distributions.}\label{Tab1}
\end{table}
\end{center}
\begin{center}
\begin{table}
\begin{center}
  \begin{tabular}{p{5.5cm}p{6.0cm}p{3.0cm}}
    \hline
    \textbf{Distribution name} &  \textbf{The CDF}  & \textbf{Domain}  \\
    \hline
    Normal distribution & $F(x)=\Phi\left(\frac{x-\mu}{\sigma}\right)$ & $x\in\mathbb{R}$ \\
    Cauchy distribution & $F(x)=\frac{1}{\pi}\tan^{-1}\left(\frac{x-\tilde{x}}{c}\right)+\frac{1}{2}$  & $x\in\mathbb{R}$ \\
    Skew-normal distribution & $F(x)=\Phi\left(\frac{x-\eta}{\varpi}\right)-2\mathcal{T}\left(\frac{x-\eta}{\varpi},\varrho\right)$& $x\in\mathbb{R}$ \\
    Gamma distribution &  $F(x)=\frac{\gamma(\kappa,x/\theta)}{\Gamma(\kappa)}$ &$x\in[0,+\infty)$\\
    \hline
  \end{tabular}
\end{center}
 \caption{The CDF of some distributions.}\label{Tab2}
\end{table}
\end{center}
\begin{figure}[h]
  \hspace{-0.0cm}\includegraphics[width=1.5\linewidth]{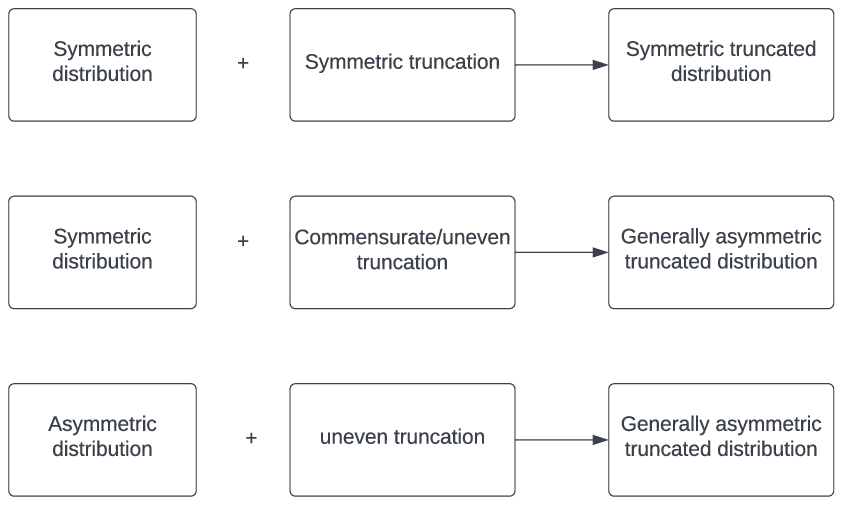}
  \vspace{-5.0cm}\caption{\footnotesize A scheme illustrates the effects of the truncations on various distributions.}\label{Fig2}
\end{figure}
By applying a truncation on a distribution defined on the real line, the resulting truncated distribution can be symmetric or asymmetric, depending on both. For instance, when the symmetric truncation abridges a symmetric distribution, the eventual truncated distribution is symmetric, while the commensurate truncation yields asymmetric truncated distribution. The possible cases of applying various truncations on symmetric/asymmetric distributions are shown in Figure \ref{Fig2}. Before commencing the mathematical definition and analysis of truncation, we mention that symmetric truncation is a particular case of commensurate truncation, so we investigate the commensurate and uneven cases in detail.\\
 First, we start with the commensurate case, i.e., when the number of deleted subintervals in the left and right parts of the real line is equal ($\tilde{M}=N$). By implementing the basic definition of the truncated probability distribution \cite{thomopoulos2018probability} and extending it to $\hat{N}=2N$ intervals, we have
\begin{equation}\label{Truncone}
  g_{X}(x)=\frac{f_{X}(x)}{\mathcal{A}},~~x\in\Omega,
\end{equation}
where
\begin{equation}\label{DDF1}
  \begin{array}{l}
    \mathcal{A}=F(b)-\mathcal{A}^{-}-\mathcal{A}^{+}-F(a), \\
    \mathcal{A}^{-}=\displaystyle{\sum_{j=1}^{N}}\left(F(\zeta_{j})-F(\vartheta_{j})\right),~~
    \mathcal{A}^{+}=\displaystyle{\sum_{j=1}^{N}}\left(F(\beta_{j})-F(\alpha_{j})\right).
  \end{array}
\end{equation}
Next, we present the associated distribution function considering these truncated intervals. Since we have disjoint closed intervals (almost everywhere) where the hidden target exists in one of them, then we have a finite sequence of distribution functions $\{G_{m}(x)\}_{m=0}^{\hat{N}+1}$ corresponding to the intervals $\{\Omega_{m}\}_{m=0}^{\hat{N}+1}$ which are listed from far left to far right, such that
\begin{equation}\label{omegatwo}
\begin{array}{l}
\{\Omega_{j}\}_{j=0}^{N-1}=\{[\zeta_{N+1-j},\vartheta_{N-j}]\}_{j=0}^{N-1},~\Omega_{N}=[\zeta_{1},0],\\\\
 \Omega_{N+1}=[0,\alpha_{1}],~\{\Omega_{N+j}\}_{j=2}^{N+1}=\{[\beta_{j-1},\alpha_{j}]\}_{j=2}^{N+1},
\end{array}
\end{equation}
and $G_{m}(x)$, $m=0,1,2,\ldots,2N+1$ are given by
\begin{equation}\label{Gfuns1}
  \begin{array}{rl}
    G_{0}(x)= & \frac{F(x)-F(a)}{\mathcal{A}},~~x\in\Omega_{0}, \\
     G_{m}(x)= & \frac{1}{\mathcal{A}}\left(F(x)-\displaystyle{\sum_{j=0}^{m-1}}(F(\zeta_{N-j})-F(\vartheta_{N-j}))-F(a)\right),~~\\
    & x\in\Omega_{m},~m=1,2,\ldots,N-1, \\
    G_{N}(x)= & \frac{F(x)-\mathcal{A}^{-}-F(a)}{\mathcal{A}},~~x\in\Omega_{N}, \\
    G_{N+1}(x)= &G_{N}(x),~~x\in\Omega_{N+1}, \\
    G_{m}(x)= &  \frac{1}{\mathcal{A}}\left(F(x)-\mathcal{A}^{-}-\displaystyle{\sum_{j=1}^{m-1-N}}
    \Big(F(\beta_{j})-F(\alpha_{j})\Big)-F(a)\right),~~\\
    &x\in\Omega_{m},~~m=N+2,~N+3,\ldots,\hat{N}+1.
  \end{array}
\end{equation}
Next, for the uneven truncation case, i.e., the number of omitting subintervals of negative and positive parts of the real line is unequal ($\tilde{M}=M$). The probability distribution function  with $\hat{N}=N+M$ is given by
\begin{equation}\label{Trunctwo}
   h_{X}(x)=\frac{f_{X}(x)}{\mathcal{A}_{\mbox{asy}}},~~x\in\tilde{\Omega},
\end{equation}
where
$$\tilde{\Omega}=\left(\bigcup_{j=1}^{M}[\zeta_{j+1},\vartheta_{j}]\right)
\cup[\zeta_{1},\alpha_{1}]\cup\left(\displaystyle{\bigcup_{\ell=1}^{N}}[\beta_{\ell},\alpha_{\ell+1}]\right),$$
$$ \begin{array}{l}
    \mathcal{A}_{\mbox{asy}}=F(b)-\mathcal{A}_{\mbox{asy}}^{-}-\mathcal{A}_{\mbox{asy}}^{+}-F(a), \\
    \mathcal{A}_{\mbox{asy}}^{-}=\displaystyle{\sum_{j=1}^{M}}\left(F(\zeta_{j})-F(\vartheta_{j})\right),~~
    \mathcal{A}_{\mbox{asy}}^{+}=\displaystyle{\sum_{j=1}^{N}}\left(F(\beta_{j})-F(\alpha_{j})\right).
  \end{array}$$
 And its cumulative distribution function is given by
\begin{equation}\label{Hfuns2}
  \begin{array}{rl}
    H_{0}(x)= & \frac{F(x)-F(a)}{\mathcal{A}},~~x\in\tilde{\Omega}_{0}, \\
     H_{m}(x)= & \frac{1}{\mathcal{A}}\left(F(x)-\displaystyle{\sum_{j=0}^{m-1}}(F(\zeta_{M-j})-F(\vartheta_{M-j}))-F(a)\right),~~\\
    & x\in\tilde{\Omega}_{m},~m=1,2,\ldots,M-1, \\
    H_{M}(x)= & \frac{F(x)-\mathcal{A}^{-}-F(a)}{\mathcal{A}},~~x\in\tilde{\Omega}_{M}, \\
    H_{M+1}(x)= & H_{M}(x),~~x\in\tilde{\Omega}_{M+1}, \\
    H_{m}(x)= &  \frac{1}{\mathcal{A}}\left(F(x)-\mathcal{A}^{-}-\displaystyle{\sum_{j=1}^{m-1-M}}
    \Big(F(\beta_{j})-F(\alpha_{j})\Big)-F(a)\right),~~\\
    &x\in\tilde{\Omega}_{m},~~m=M+2,~M+3,\ldots,\hat{N}+1.
  \end{array}
\end{equation}
Moreover, for distributions defined on the interval $[0,+\infty)$, implementing $N$-deleted subintervals, consequently, the corresponding truncated probability distribution function is given by
\begin{equation}\label{Trunctwo}
   u_{X}(x)=\frac{f_{X}(x)}{\mathcal{A}_{r}},~~x\in\hat{\Omega},
\end{equation}
such that
$$\hat{\Omega}=[0,\alpha_{1}]\cup\left(\displaystyle{\bigcup_{\ell=1}^{N}}[\beta_{\ell},\alpha_{\ell+1}]\right),$$
$$
    \mathcal{A}_{r}=F(b)-\mathcal{A}_{r}^{+},~~~\mathcal{A}_{r}^{+}=\sum_{j=1}^{N}\left(F(\beta_{j})-F(\alpha_{j})\right),
$$
  and its accumulated distribution function is
  \begin{equation}\label{Ufuns3}
  \begin{array}{rl}
    U_{0}(x)= & \frac{F(x)}{\mathcal{A}_{r}} ,~~x\in\hat{\Omega}_{0}=[0,\alpha_{1}], \\
    U_{m}(x)= &  \frac{1}{\mathcal{A}_{r}}\left(F(x)-\displaystyle{\sum_{j=1}^{m}}
    \Big(F(\beta_{j})-F(\alpha_{j})\Big)\right),~~\\
    &x\in\hat{\Omega}_{m},~~m=1,~2,\ldots,N.
  \end{array}
\end{equation}
\begin{figure}[h]
  \hspace{-0.0cm}\includegraphics[width=1.3\linewidth]{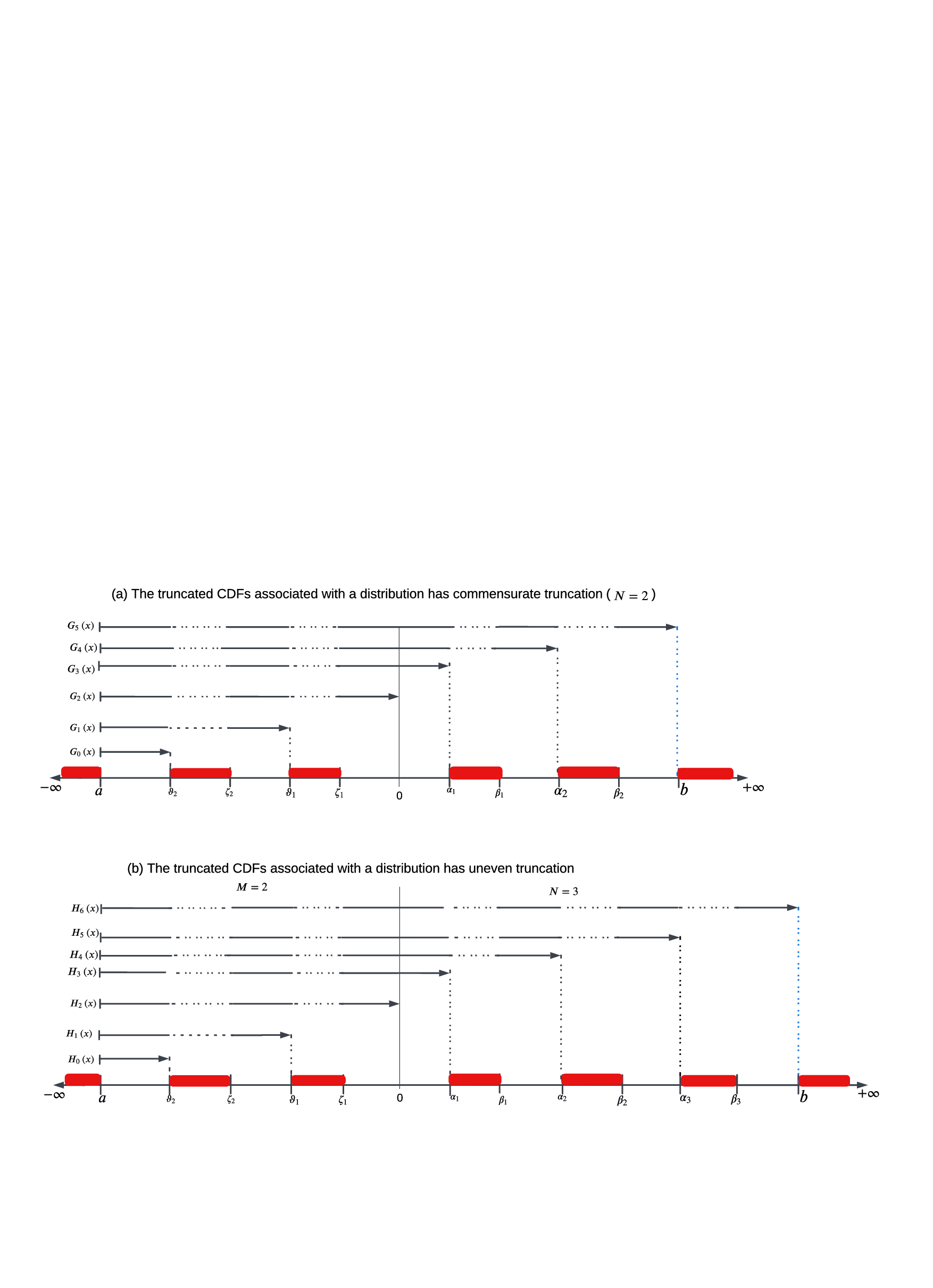}
  \vspace{-1.5cm}\caption{\footnotesize The domain of truncated CDF under several type of truncations.}\label{Fig3}
\end{figure}
\section{Optimization of the searching parameters}
One of the crucial factors in the searching process is estimating the elapsed time to find the lost target and optimizing it. Here, we have two searchers move with the same uniform velocity $v_{1}$ within the searching regions, and with $v_{2}$ in the deleted regions, the lost target is guaranteed to do not exist. Consequently, the elapsed time to sweep a searching region is the length of this interval divided by $v_{1}$, and in a deleted interval, it takes time equal to the measure of this interval divided by $v_{2}$. We establish a finite sequence of elapsed time $\{\tau_{m}\}_{m=0}^{\hat{N}+1}$, considering the hidden target lies in one of the intervals $\Omega_{j}$, $j=0,1,\ldots \hat{N}+1$, ($\hat{N}=2N,~N+M$ for the commensurate and uneven cases) the communication between the two searchers, and the searching process is stopped, once the lost target is detected. The elapsed time $\tau_{m}$ for searches sweep commensurate truncated domain is given by
\begin{equation}\label{timeeq1}
\tau_{m}=\left\{\begin{array}{rl}
                  \frac{1}{v_{1}}\displaystyle{\sum_{\ell=1}^{N+1-m}}(\vartheta_{\ell-1}-\zeta_{\ell})+\frac{\delta_{N+1-m,1}}{v_{2}}
                  \displaystyle{\sum_{\ell=1}^{N-m}}(\zeta_{\ell}-\vartheta_{\ell}), & 0\leq m\leq N, \\
                  \frac{1}{v_{1}}\displaystyle{\sum_{\ell=1}^{m-N}}(\alpha_{\ell}-\beta_{\ell-1})+\frac{\delta_{m-N,1}}{v_{2}}
                  \displaystyle{\sum_{\ell=1}^{m-N-1}}(\beta_{\ell}-\alpha_{\ell}), & N< m\leq \hat{N}+1,
                \end{array}\right.
\end{equation}
where $\delta_{k,1}$ is defined by
$$\delta_{k,1}=\left\{\begin{array}{rl}
                        0, & k\leq1, \\
                        1, & k>1,
                      \end{array}\right.
$$
and for the searchers scanning an uneven truncated domain  is given by
\begin{equation}\label{timeeq2}
\tau_{m}=\left\{\begin{array}{rl}
                  \frac{1}{v_{1}}\displaystyle{\sum_{\ell=1}^{M+1-m}}(\vartheta_{\ell-1}-\zeta_{\ell})+\frac{\delta_{M+1-m,1}}{v_{2}}
                  \displaystyle{\sum_{\ell=1}^{M-m}}(\zeta_{\ell}-\vartheta_{\ell}), & 0\leq m\leq M, \\
                  \frac{1}{v_{1}}\displaystyle{\sum_{\ell=1}^{m-M}}(\alpha_{\ell}-\beta_{\ell-1})+\frac{\delta_{m-M,1}}{v_{2}}
                  \displaystyle{\sum_{\ell=1}^{m-M-1}}(\beta_{\ell}-\alpha_{\ell}), & M< m\leq \hat{N}+1.
                \end{array}\right.
\end{equation}
For a searcher $S$ with probability distribution function defined over $[0,+\infty)$, the estimated time $\{\tau_{m}\}_{m=1}^{N+1}$ to find a lost target is calculated by
\begin{equation}\label{timeeq3}
\tau_{m}=\frac{1}{v_{1}}\displaystyle{\sum_{\ell=1}^{m}}(\alpha_{\ell}-\beta_{\ell-1})+\frac{\delta_{m,1}}{v_{2}}
                  \displaystyle{\sum_{\ell=1}^{m-1}}(\beta_{\ell}-\alpha_{\ell}), ~~~ 1\leq m\leq N+1.
\end{equation}
Consequently, the expected elapsed time $\mathbb{E}(\tau_{m})$ is established based on equations \eqref{Gfuns1} and \eqref{timeeq1} for distributions undergo commensurate truncation such that
\begin{equation}\label{expcttmeq1}
 \mathbb{E}(\tau_{m})=\left\{\begin{array}{rl}
                  \frac{1}{v_{1}}\displaystyle{\sum_{\ell=1}^{N+1-m}}(\vartheta_{\ell-1}-\zeta_{\ell})
                  \mathbb{P}_{N+1-\ell}\Big([\zeta_{\ell},\vartheta_{\ell-1}]\Big)+\frac{\delta_{N+1-m,1}}{v_{2}}
                  \displaystyle{\sum_{\ell=1}^{N-m}}(\zeta_{\ell}-\vartheta_{\ell}), & 0\leq m\leq N, \\
                  \frac{1}{v_{1}}\displaystyle{\sum_{\ell=1}^{m-N}}(\alpha_{\ell}-\beta_{\ell-1})
                  \mathbb{P}_{N+\ell}\Big([\beta_{\ell-1},\alpha_{\ell}]\Big)+\frac{\delta_{m-N,1}}{v_{2}}
                  \displaystyle{\sum_{\ell=1}^{m-N-1}}(\beta_{\ell}-\alpha_{\ell}), & N< m\leq 2N+1,
                \end{array}\right.
\end{equation}
where
\begin{equation}\label{ProbG}
\begin{array}{rl}
    \mathbb{P}_{N+1-\ell}\Big([\zeta_{\ell},\vartheta_{\ell-1}]\Big)=G_{N+1-\ell}(\vartheta_{\ell-1})-G_{N+1-\ell}(\zeta_{\ell}), & 0\leq m\leq N, \\
    \mathbb{P}_{N+\ell}\Big([\beta_{\ell-1},\alpha_{\ell}]\Big)=G_{N+\ell}(\alpha_{\ell})-G_{N+\ell}(\beta_{\ell-1}), & N< m\leq \hat{N}+1.
  \end{array}
\end{equation}
Next based on equations \eqref{Hfuns2} and \eqref{timeeq2},  the expected time to detect a hidden target by two searchers sweep uneven truncated domain is given by
\begin{equation}\label{expcttmeq2}
 \mathbb{E}(\tau_{m})=\left\{\begin{array}{rl}
                  \frac{1}{v_{1}}\displaystyle{\sum_{\ell=1}^{M+1-m}}(\vartheta_{\ell-1}-\zeta_{\ell})
                  \mathbb{P}_{M+1-\ell}\Big([\zeta_{\ell},\vartheta_{\ell-1}]\Big)+\frac{\delta_{M+1-m,1}}{v_{2}}
                  \displaystyle{\sum_{\ell=1}^{M-m}}(\zeta_{\ell}-\vartheta_{\ell}), & 0\leq m\leq M, \\
                  \frac{1}{v_{1}}\displaystyle{\sum_{\ell=1}^{m-M}}(\alpha_{\ell}-\beta_{\ell-1})
                  \mathbb{P}_{M+\ell}\Big([\beta_{\ell-1},\alpha_{\ell}]\Big)+\frac{\delta_{m-M,1}}{v_{2}}
                  \displaystyle{\sum_{\ell=1}^{m-M-1}}(\beta_{\ell}-\alpha_{\ell}), & M< m\leq \hat{N}+1,
                \end{array}\right.
\end{equation}
\begin{equation}\label{ProbH}
\begin{array}{rl}
    \mathbb{P}_{M+1-\ell}\Big([\zeta_{\ell},\vartheta_{\ell-1}]\Big)=H_{M+1-\ell}(\vartheta_{\ell-1})-H_{M+1-\ell}(\zeta_{\ell}), & 0\leq m\leq M, \\
    \mathbb{P}_{M+\ell}\Big([\beta_{\ell-1},\alpha_{\ell}]\Big)=H_{M+\ell}(\alpha_{\ell})-H_{M+\ell}(\beta_{\ell-1}), & M< m\leq \hat{N}+1.
  \end{array}
\end{equation}
Finally, for one searcher on a half domain $[0,\infty)$ with $N$-truncated intervals, its expected elapsed time to find the hidden target is obtained by incorporating equations \eqref{Ufuns3} and \eqref{timeeq3} such that
\begin{equation}\label{expcttmeq3}
\mathbb{E}(\tau_{m})=\frac{1}{v_{1}}\displaystyle{\sum_{\ell=1}^{m}}(\alpha_{\ell}-\beta_{\ell-1})
\mathbb{P}_{\ell-1}\Big([\beta_{\ell-1},\alpha_{\ell}]\Big)+\frac{\delta_{m,1}}{v_{2}}
                  \displaystyle{\sum_{\ell=1}^{m-1}}(\beta_{\ell}-\alpha_{\ell}), ~~~ 1\leq m\leq N+1,
\end{equation}
where
\begin{equation}\label{ProbU}
\mathbb{P}_{\ell-1}\Big([\beta_{\ell-1},\alpha_{\ell}]\Big)=U_{\ell-1}(\alpha_{\ell})-U_{\ell-1}(\beta_{\ell-1}).
\end{equation}
Next, we will estimate the optimum expected elapsed time for each case. Here, we consider the expected elapsed time function $\mathbb{E}(\tau_{m})$ as the objective function subjected to the following constraints
\begin{equation}\label{cnstobj1}
  \begin{array}{l}
    0<-\zeta_{i}<-\vartheta_{i},~~~0<\alpha_{i}<\beta_{i},~~i=1,2,\ldots,N, \mbox{ for commensurate case,} \\
     0<-\zeta_{j}<-\vartheta_{j},~j=1,2,\ldots,M,~~0<\alpha_{i}<\beta_{i},~~i=1,2,\ldots,N, \mbox{ for uneven case,}\\
    0<\alpha_{i}<\beta_{i},~~i=1,2,\ldots,N, \mbox{ for } N\mbox{-deleted subintervals on }  [0,\infty).
  \end{array}
\end{equation}
On the one hand, to obtain the optimum value of $\mathbb{E}(\tau_{m})$ considering all the truncations in the negative part of the real line, we investigate $\mathbb{E}(\tau_{0})$, i.e., $m=0$ into equations \eqref{expcttmeq1}, and \eqref{expcttmeq2}. On the other hand, to have its optimum value taking into account all the truncations in the positive part of the real line, $\mathbb{E}(\tau_{\breve{N}})$ has been studied, where $\breve{N}=2N+1$, $N+M+1$, and $N+1$ into equations \eqref{expcttmeq1}, \eqref{expcttmeq2}, and \eqref{expcttmeq3}, respectively. Consequently, the nonlinear programming is to optimize the objective functions $\mathbb{E}(\tau_{0})$, and  $\mathbb{E}(\tau_{\breve{N}})$ subjected to constraints \eqref{cnstobj1}. The optimization techniques can be classified into two main categories; unconditional and conditional optimization \cite{fiacco1990nonlinear}. Here, we are interested in investigating the objective functions in view of unconditional optimization. To achieve this aim, we incorporate conditions \eqref{cnstobj1} with equations \eqref{expcttmeq1}, \eqref{expcttmeq2}, and \eqref{expcttmeq3}, using the artificial parameters $\left\{\frac{\epsilon_{\ell}}{v_{2}^{2}}\right\}_{\ell=1}^{\tilde{M}}$, $\left\{\frac{\hat{\epsilon}_{\ell}}{v_{1}^{2}}\right\}_{\ell=1}^{\tilde{M}}$, and $\left\{\frac{\tilde{\epsilon}_{\ell}}{v_{1}^{2}}\right\}_{\ell=1}^{\tilde{M}}$, with the objective functions $\mathbb{E}(\tau_{0})$ of symmetric and asymmetric cases, and the parameters $\left\{\frac{\varepsilon_{\ell}}{v_{2}^{2}}\right\}_{\ell=1}^{N}$, $\left\{\frac{\hat{\varepsilon}_{\ell}}{v_{1}^{2}}\right\}_{\ell=1}^{N}$, and $\left\{\frac{\tilde{\varepsilon}_{\ell}}{v_{1}^{2}}\right\}_{\ell=1}^{N}$ with $\mathbb{E}(\tau_{\breve{N}})$ of the three cases. Consequently, we have the following modified objective functions
\begin{equation}\label{mdfEsym1}
\begin{array}{rl}
  \hat{\mathbb{E}}(\tau_{0})= & \frac{1}{v_{1}}\displaystyle{\sum_{\ell=1}^{\tilde{M}+1}}(\vartheta_{\ell-1}-\zeta_{\ell})
                  \mathbb{P}_{\tilde{M}+1-\ell}\Big([\zeta_{\ell},\vartheta_{\ell-1}]\Big)+\frac{1}{v_{2}}
                  \displaystyle{\sum_{\ell=1}^{\tilde{M}}}\Big(\zeta_{\ell}-\vartheta_{\ell}+
                  \frac{\epsilon_{\ell}}{v_{2}}(\zeta_{\ell}-\vartheta_{\ell})^{2}\Big)+ \\
&\frac{1}{v_{1}^{2}}\displaystyle{\sum_{\ell=1}^{\tilde{M}}}(\hat{\epsilon}_{\ell}\zeta^{2}_{\ell}+
\tilde{\epsilon}_{\ell}\vartheta^{2}_{\ell}),~~\tilde{M}=N,~M,\\
\end{array}
\end{equation}
\begin{equation}\label{mdfEsym2}
\begin{array}{rl}
  \hat{\mathbb{E}}(\tau_{\breve{N}})=&\frac{1}{v_{1}}\displaystyle{\sum_{\ell=1}^{N+1}}(\alpha_{\ell}-
  \beta_{\ell-1})\mathbb{P}_{\tilde{M}+\ell}\Big([\beta_{\ell-1},\alpha_{\ell}]\Big)+\frac{1}{v_{2}} \displaystyle{\sum_{\ell=1}^{N}}\Big(\beta_{\ell}-\alpha_{\ell}+\frac{\varepsilon_{\ell}}{v_{2}}(\beta_{\ell}-\alpha_{\ell})^{2}\Big)+ \\
&\frac{1}{v_{1}^{2}}\displaystyle{\sum_{\ell=1}^{N}}(\hat{\varepsilon}_{\ell}\beta^{2}_{\ell}+
\tilde{\varepsilon}_{\ell}\alpha^{2}_{\ell}),~~\breve{N}=2N+1,~N+M+1,\\
\end{array}
\end{equation}
and
\begin{equation}\label{mdfEsym3}
\begin{array}{rl}
  \hat{\mathbb{E}}(\tau_{N+1})=&\frac{1}{v_{1}}\displaystyle{\sum_{\ell=1}^{N+1}}(\alpha_{\ell}-
  \beta_{\ell-1})\mathbb{P}_{\ell-1}\Big([\beta_{\ell-1},\alpha_{\ell}]\Big)+\frac{1}{v_{2}} \displaystyle{\sum_{\ell=1}^{N}}\Big(\beta_{\ell}-\alpha_{\ell}+\frac{\varepsilon_{\ell}}{v_{2}}(\beta_{\ell}-\alpha_{\ell})^{2}\Big)+ \\
&\frac{1}{v_{1}^{2}}\displaystyle{\sum_{\ell=1}^{N}}(\hat{\varepsilon}_{\ell}\beta^{2}_{\ell}+\tilde{\varepsilon}_{\ell}\alpha^{2}_{\ell}).\\
\end{array}
\end{equation}
Here all the parameters $\epsilon_{\ell},~\hat{\epsilon}_{\ell},~\tilde{\epsilon}_{\ell},$ $\varepsilon_{\ell},~\hat{\varepsilon}_{\ell},~\tilde{\varepsilon}_{\ell},$ are greater than $0$ and less than $1$, while the division by $v_{1}^{2}$, and $v_{2}^{2}$ are considered to preserve the dimensionality of the problem. Now we have unconstrained objective functions \eqref{mdfEsym1}-\eqref{mdfEsym3} with $2\tilde{M}$-variables that are minimized using Newton's method in optimization. To discuss the proposed optimization technique that is used to solve these equations altogether, we use the generic function $Q(\mathbf{x})$ with the argument vector $\mathbf{x}=(x_{1},x_{2},\ldots,x_{2\tilde{M}})^{T}$ ($(\ast)^{T}$ denotes the transpose operation), such that
$$
\begin{array}{rccl}
Q(\mathbf{x})=\tilde{\mathbb{E}}(\tau_{0}), & \left\{x_{i}=\vartheta_{i}\right\}_{i=1}^{N}, & \left\{x_{N+i}=\zeta_{i}\right\}_{i=1}^{N}, & \mbox{commensurate case}, \\
Q(\mathbf{x})=\tilde{\mathbb{E}}(\tau_{0}), & \left\{x_{i}=\vartheta_{i}\right\}_{i=1}^{M}, & \left\{x_{M+i}=\zeta_{i}\right\}_{i=1}^{M}, & \mbox{uneven case}, \\
Q(\mathbf{x})=\hat{\mathbb{E}}(\tau_{N+1}), & \left\{x_{i}=\alpha_{i}\right\}_{i=1}^{N}, & \left\{x_{N+i}=\beta_{i}\right\}_{i=1}^{N}, & \mbox{positive part}.
\end{array}
$$
The following  algorithm \ref{Algthm1} is applied
\begin{alg}\label{Algthm1}
\textbf{(Optimization Procedures)}\\
\begin{enumerate}
  \item Initial guess $\mathbf{x}^{(0)}$
  \item Calculate the gradient $\nabla Q=\left(\frac{\partial Q}{\partial x_{1}},~\frac{\partial Q}{\partial x_{2}},\ldots,\frac{\partial Q}{\partial x_{2\tilde{M}}}\right)^{T}$
  \item Calculate the Hessian matrix $\mathcal{H}=\left(\frac{\partial^{2}Q}{\partial x_{i}\partial x_{j}}\right),~i,j=1,2,\ldots,2\tilde{M}$
  \item For $k$ from 1 to $k^{\ast}$ do ($k^{\ast}$ represents the number of iteration where the convergence holds)
  \item Solve $\mathcal{H}(\mathbf{x}^{(k-1)})\mathbf{s}^{(k)}=-\nabla Q(\mathbf{x}^{(k-1)})$
  \item $\mathbf{x}^{(k)}=\mathbf{x}^{(k-1)}+\mathbf{s}^{(k)}$
  \item End for.
\end{enumerate}
\end{alg}
\large\section{Computational examples}
\begin{figure}[h]
  \vspace{0.0cm}\includegraphics[width=1.0\linewidth]{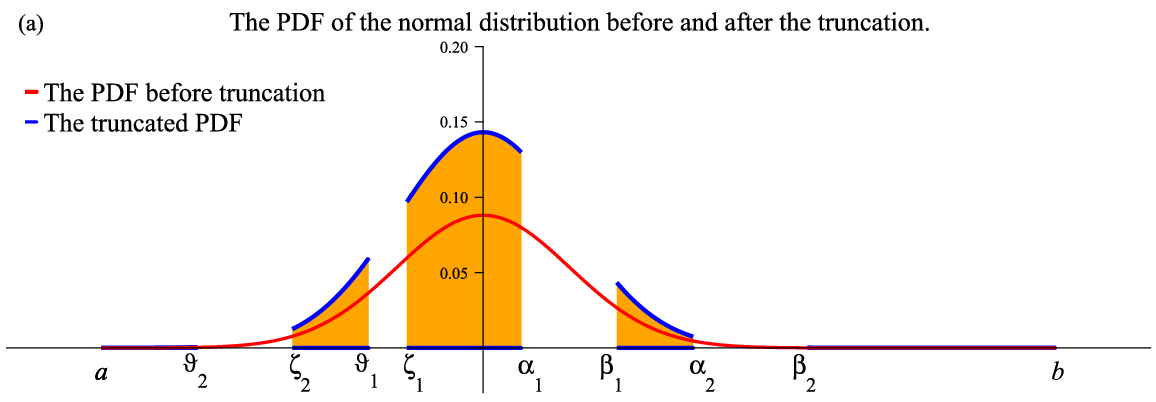}\\
  \vspace{0.0cm}\includegraphics[width=1.0\linewidth]{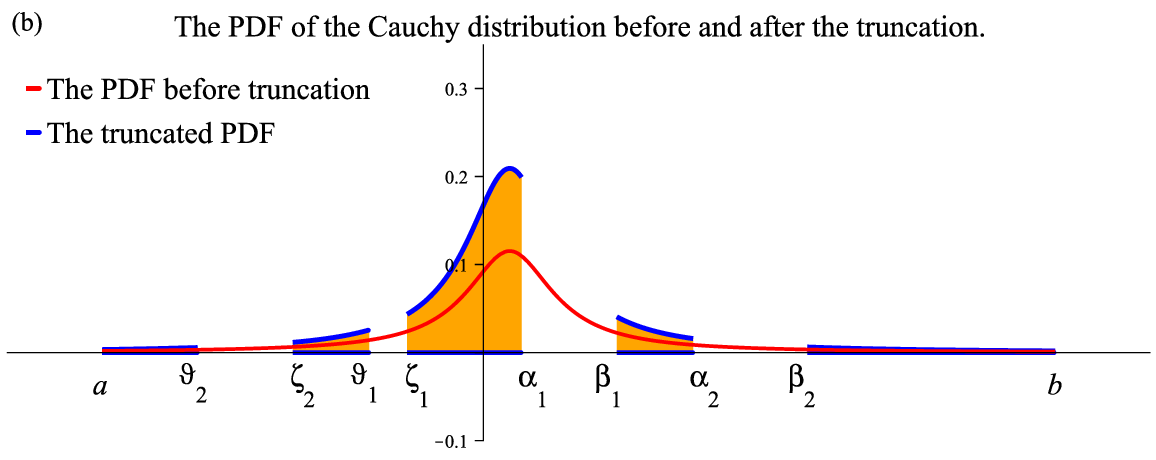}\\
   \vspace{0.0cm}\includegraphics[width=1.0\linewidth]{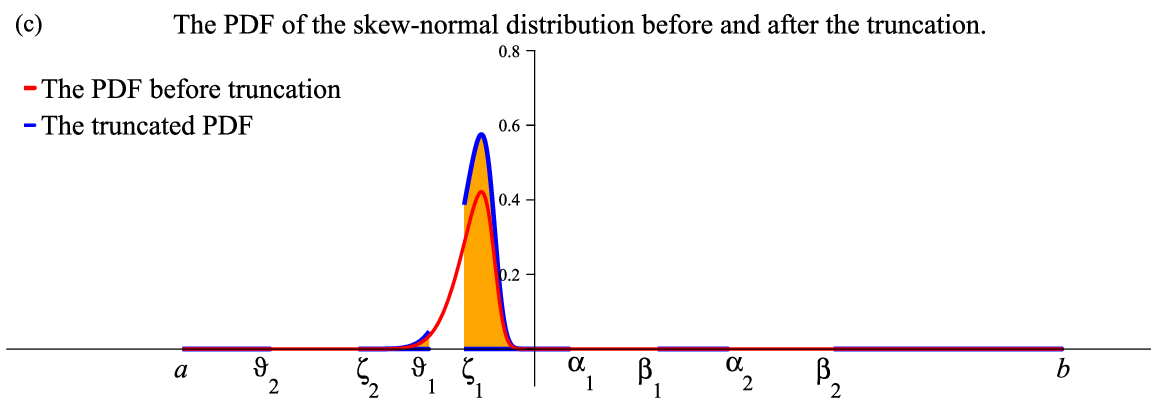}\\
  \vspace{0.0cm}\caption{\footnotesize The PDFs of the first three distributions listed in Table \ref{Tab1} and their corresponding truncated PDFs when commensurate truncation is applied with $N=2$, and the value of endpoints are $a=-20$, $b=30$.}\label{Fig4}
\end{figure}

\begin{figure}[h]
  \vspace{0.0cm}\includegraphics[width=1.0\linewidth]{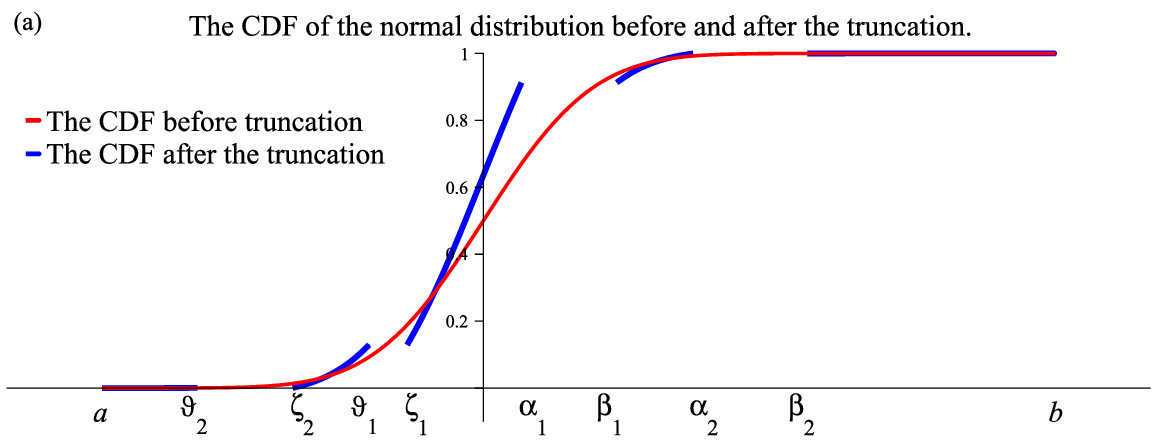}\\
  \vspace{0.0cm}\includegraphics[width=1.0\linewidth]{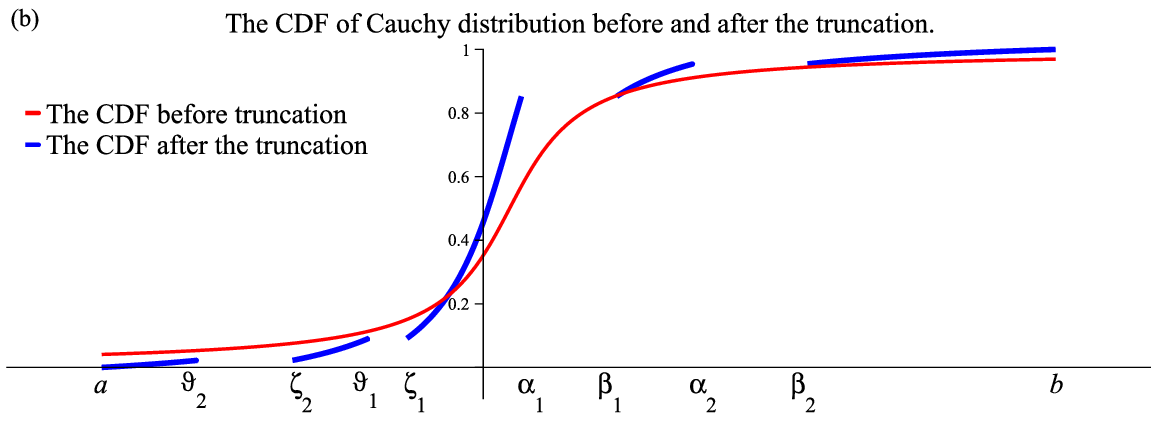}\\
   \vspace{0.0cm}\includegraphics[width=1.0\linewidth]{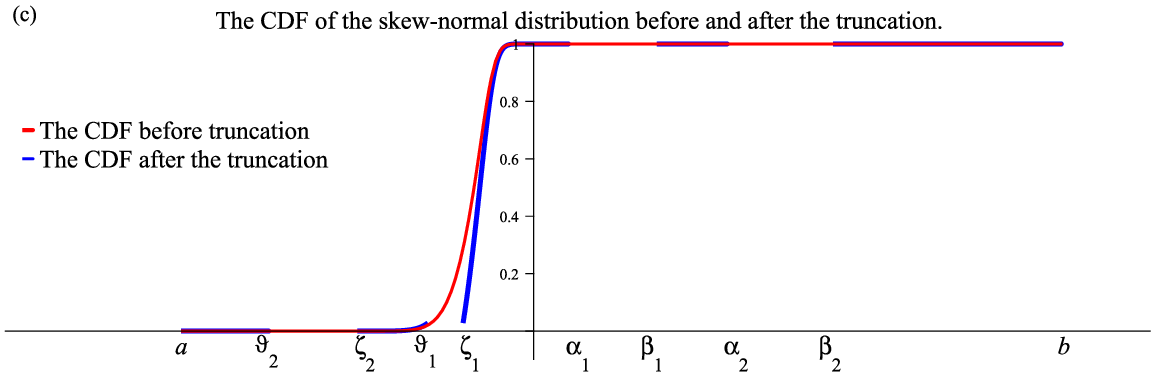}\\
  \vspace{0.0cm}\caption{\footnotesize The associated CDFs given in the first three distributions Tab \ref{Tab2}, and their related truncated CDFs under commensurate case with $N=2$, $a=-20$, $b=30$.}\label{Fig5}
\end{figure}

\begin{figure}[h]
  \vspace{0.0cm}\includegraphics[width=1.0\linewidth]{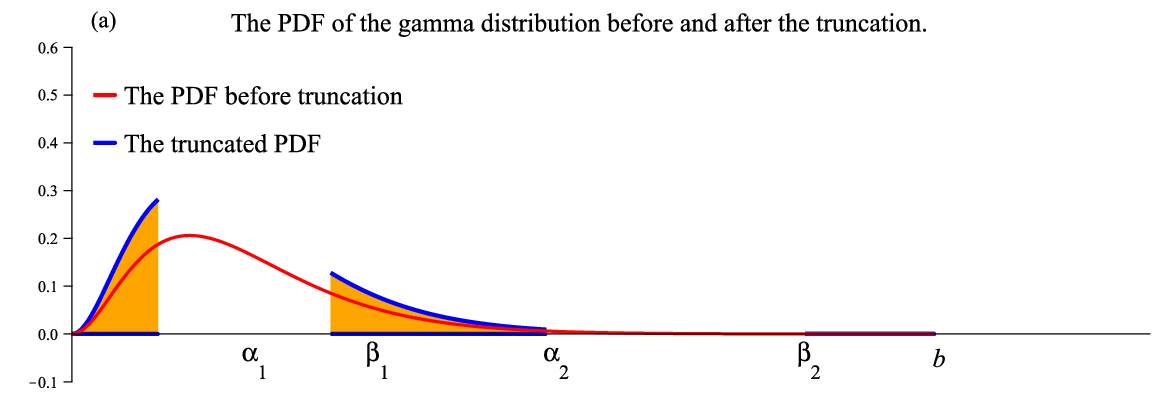}\\
  \vspace{0.0cm}\includegraphics[width=1.0\linewidth]{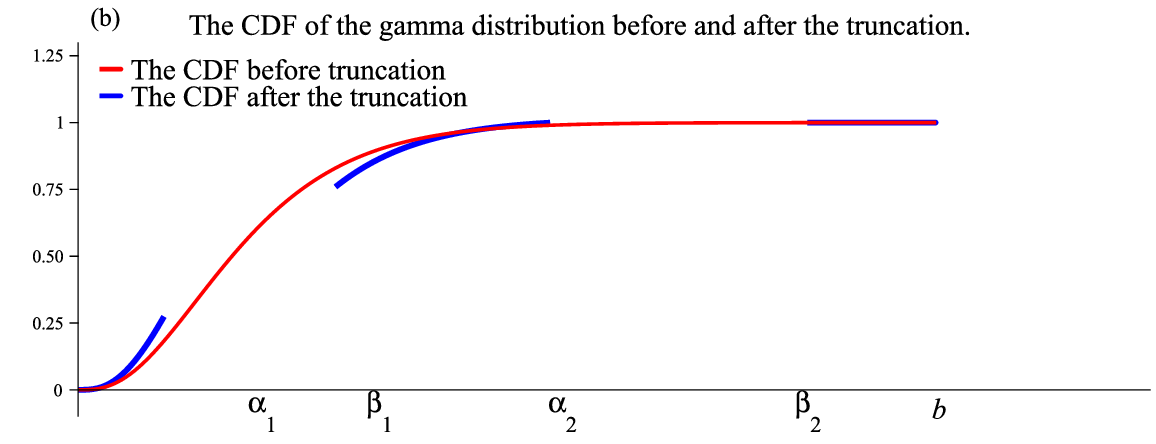}\\
  \vspace{0.0cm}\caption{\footnotesize The PDF and CDF  of he gamma distribution before and after truncation  by two deleted intervals.}\label{Fig6}
\end{figure}

\begin{figure}[h]
  \vspace{0.0cm}\includegraphics[width=1.0\linewidth]{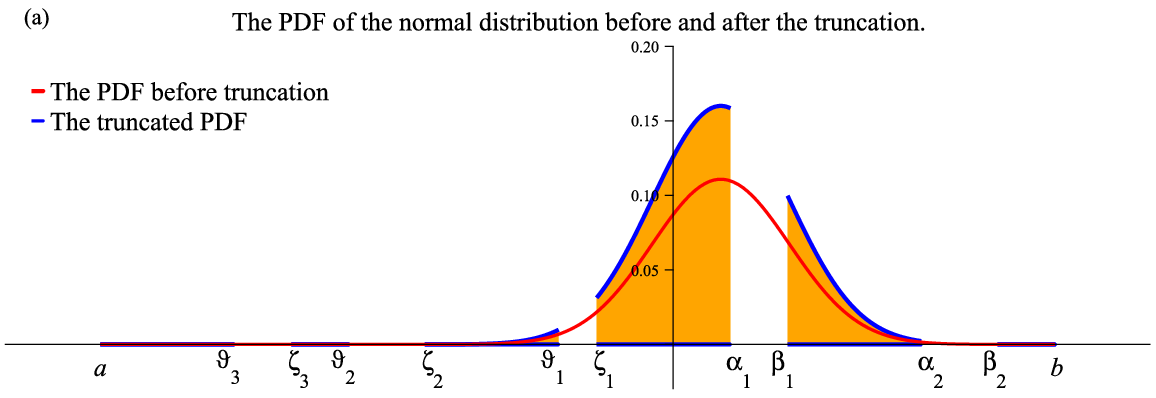}\\
  \vspace{0.0cm}\includegraphics[width=1.0\linewidth]{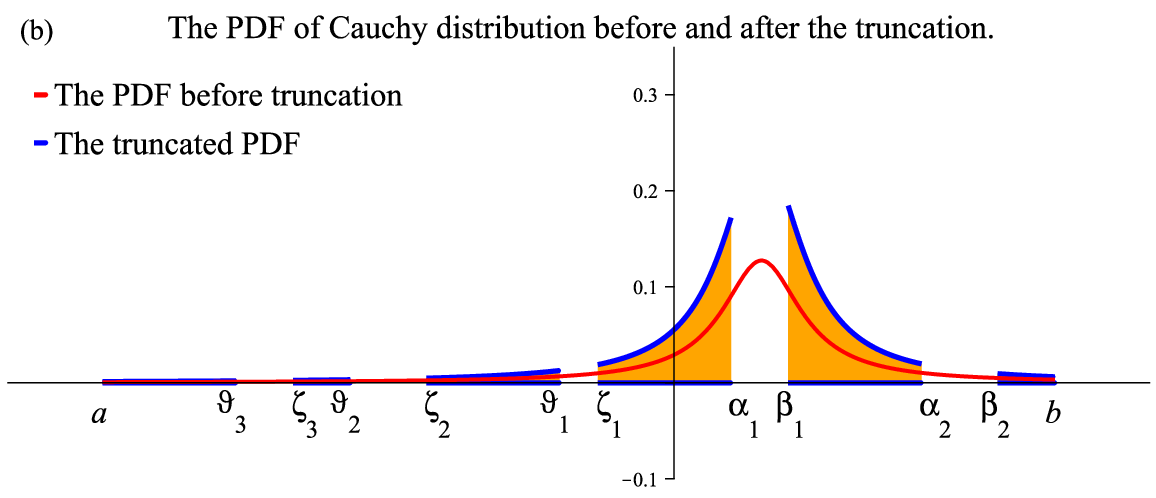}\\
   \vspace{0.0cm}\includegraphics[width=1.0\linewidth]{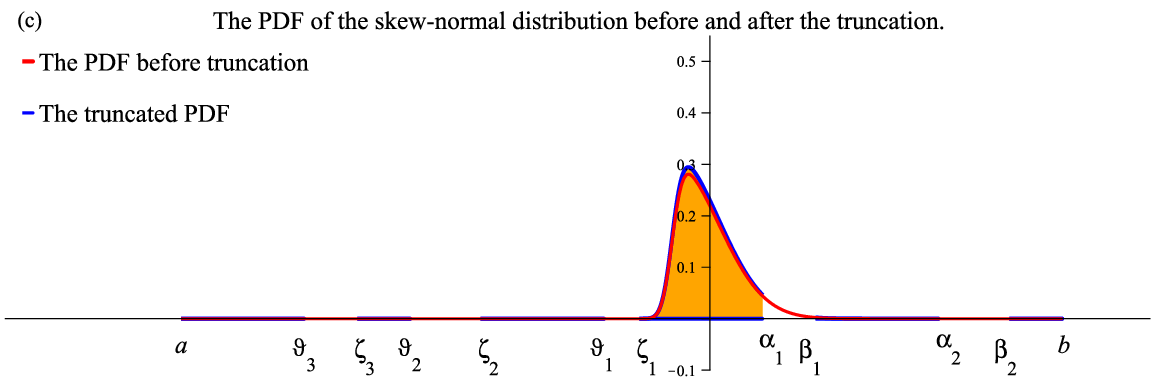}\\
  \vspace{0.0cm}\caption{\footnotesize The PDFs of the first three distributions listed in Table \ref{Tab1} and their corresponding truncated PDFs when uneven truncation is applied with $M=3$, $N=2$, and the value of endpoints are $a=-30$, $b=20$.}\label{Fig7}
\end{figure}

\begin{figure}[h]
  \vspace{0.0cm}\includegraphics[width=1.0\linewidth]{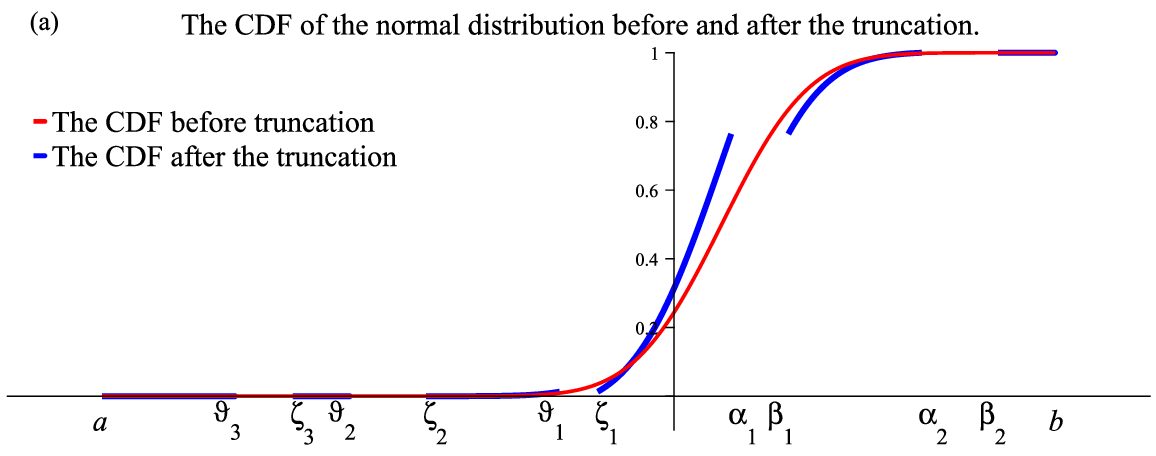}\\
  \vspace{0.0cm}\includegraphics[width=1.0\linewidth]{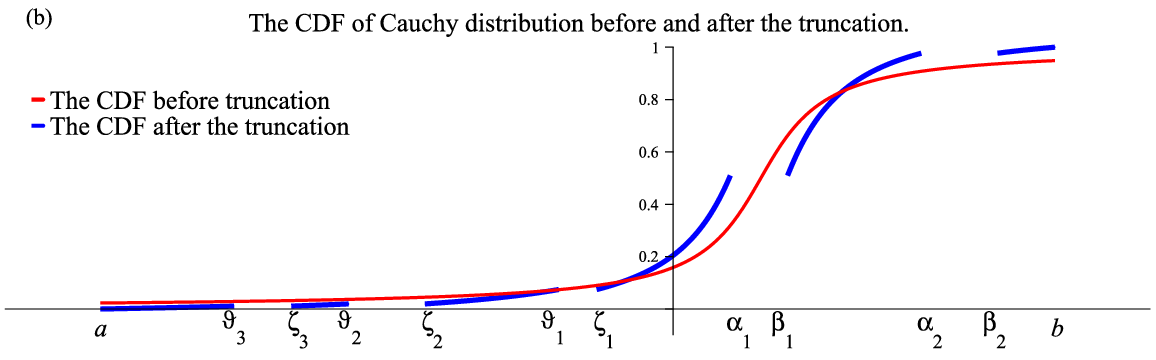}\\
   \vspace{0.0cm}\includegraphics[width=1.0\linewidth]{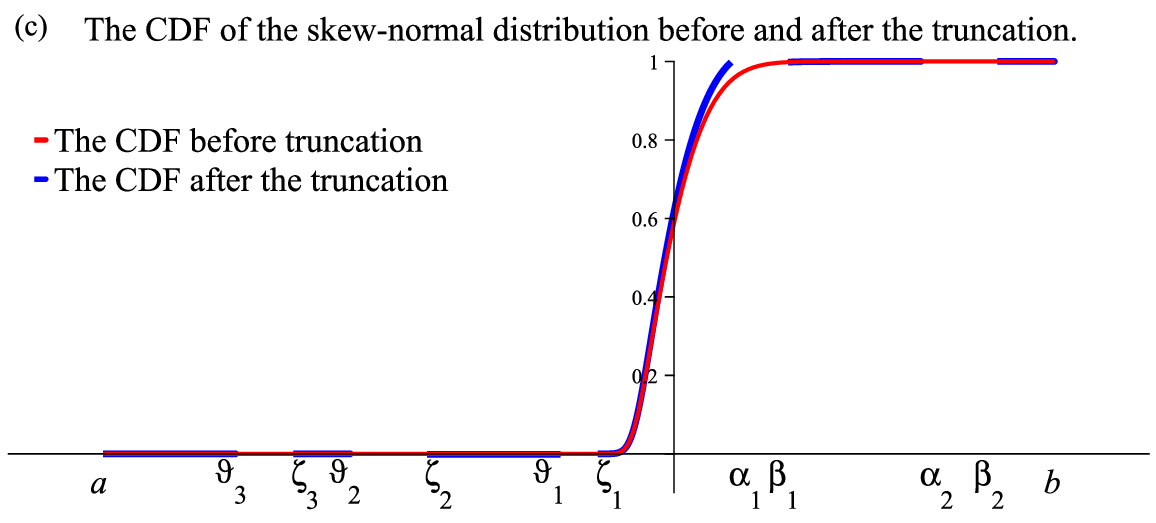}\\
  \vspace{0.0cm}\caption{\footnotesize The associated CDFs given in the first three distributions Tab \ref{Tab2}, and their related truncated CDFs  when the uneven truncation is applied with $M=3$, $N=2$, and the value of endpoints are $a=-30$, $b=20$.}\label{Fig8}
\end{figure}

\begin{figure}[h]
  \vspace{0.0cm}\includegraphics[width=1.0\linewidth]{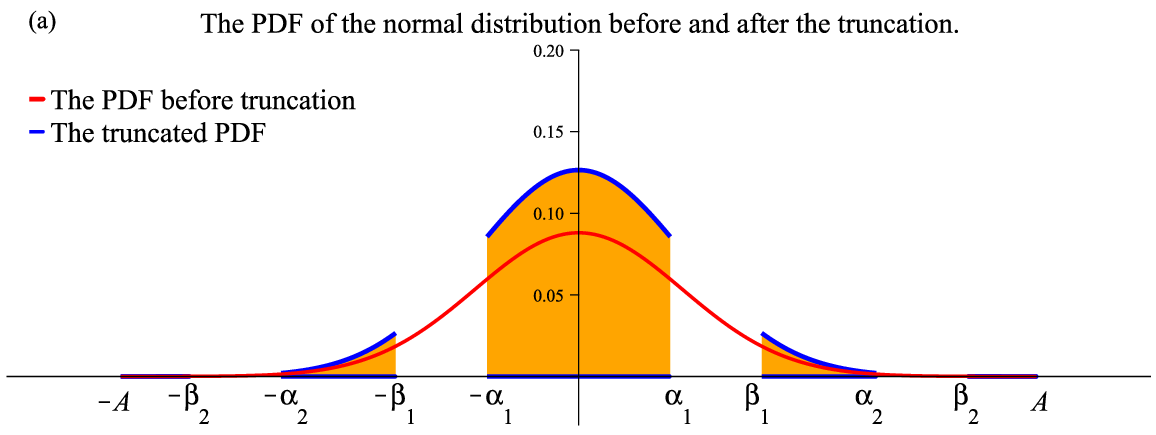}\\
  \vspace{0.0cm}\includegraphics[width=1.0\linewidth]{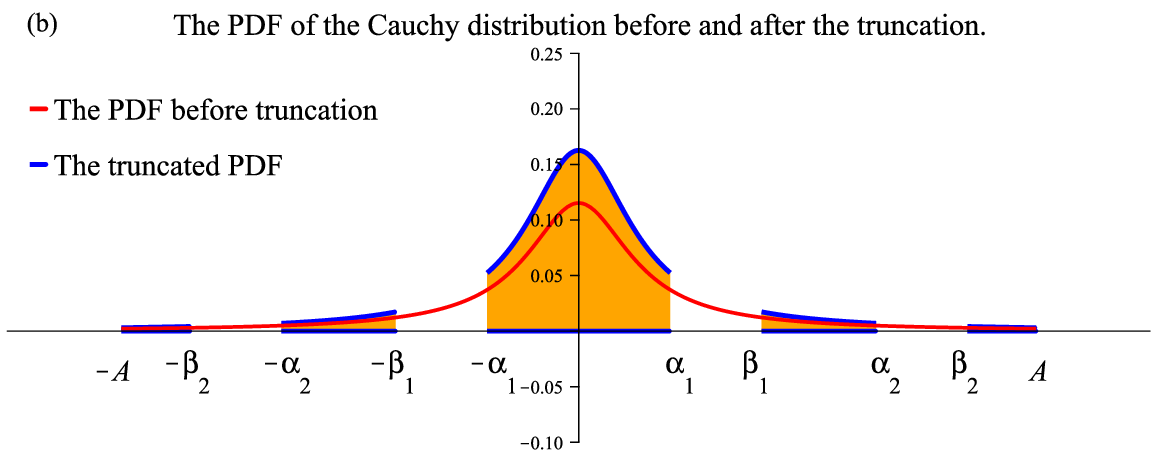}\\
  \vspace{0.0cm}\caption{\footnotesize The PDFs of the symmetric distributions listed in Table \ref{Tab1} and their corresponding truncated PDFs when the symmetric truncation is applied with $N=2$, and the value of endpoints are $b=-a=A=20$.}\label{Fig9}
\end{figure}

\begin{figure}[!htb]
   \begin{minipage}{0.5\textwidth}
     \centering
     \includegraphics[width=1.0\linewidth]{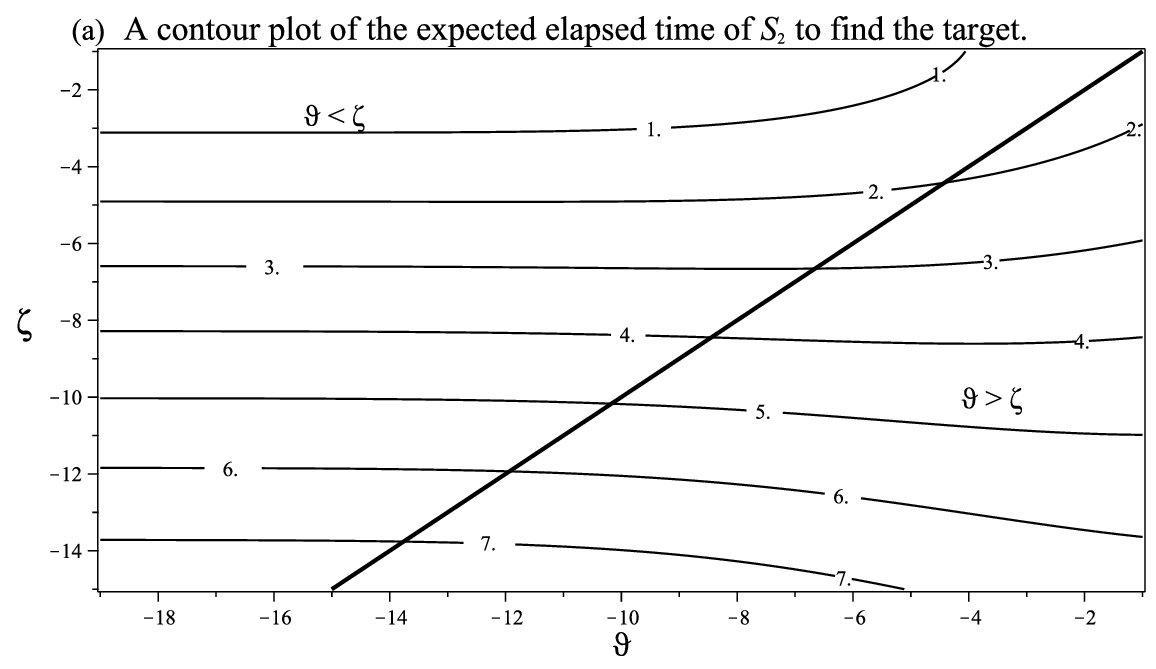}
   \end{minipage}\hfill
   \begin{minipage}{0.5\textwidth}
     \centering
  \hspace{0.5cm}   \includegraphics[width=1.0\linewidth]{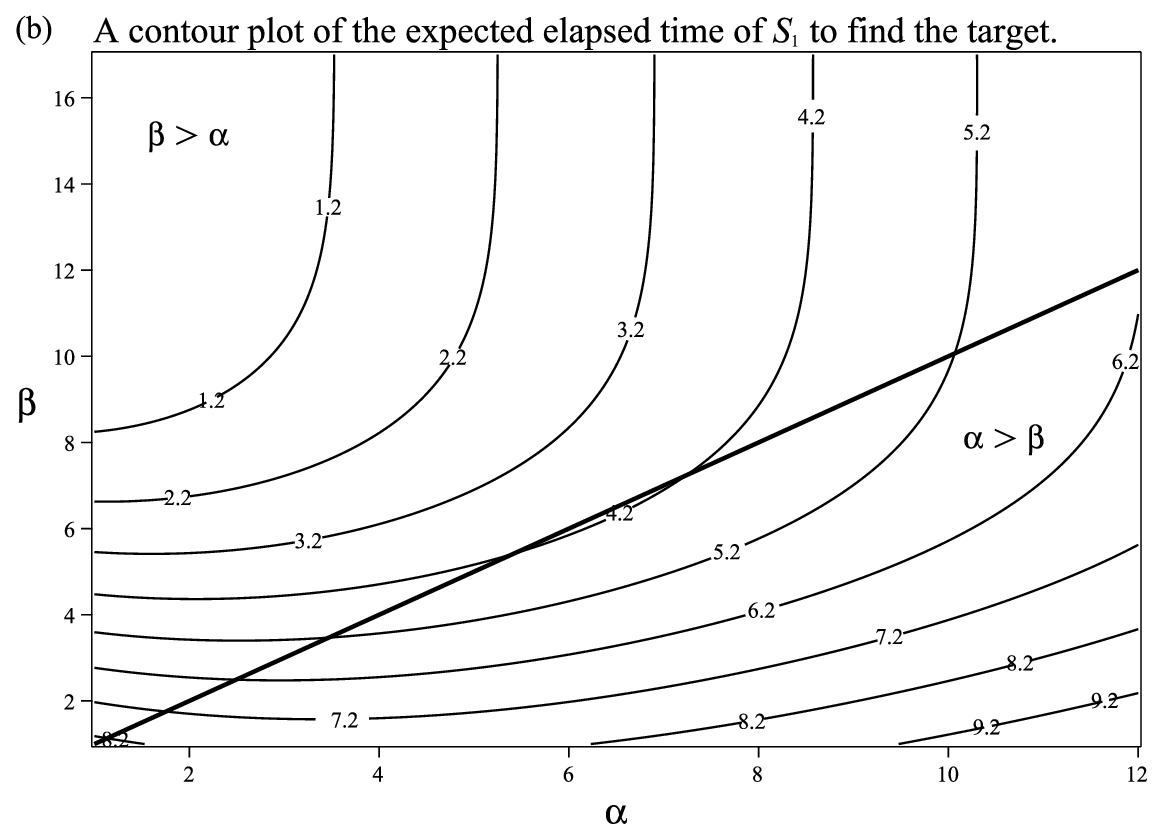}
   \end{minipage}
     \caption{A contour plot of the expected time searchers $S_{1}$, and $S_{2}$ of a target follows normal distribution.}\label{Fig10}
\end{figure}

\begin{figure}[h]
   \vspace{0.0cm}\includegraphics[width=1.0\linewidth]{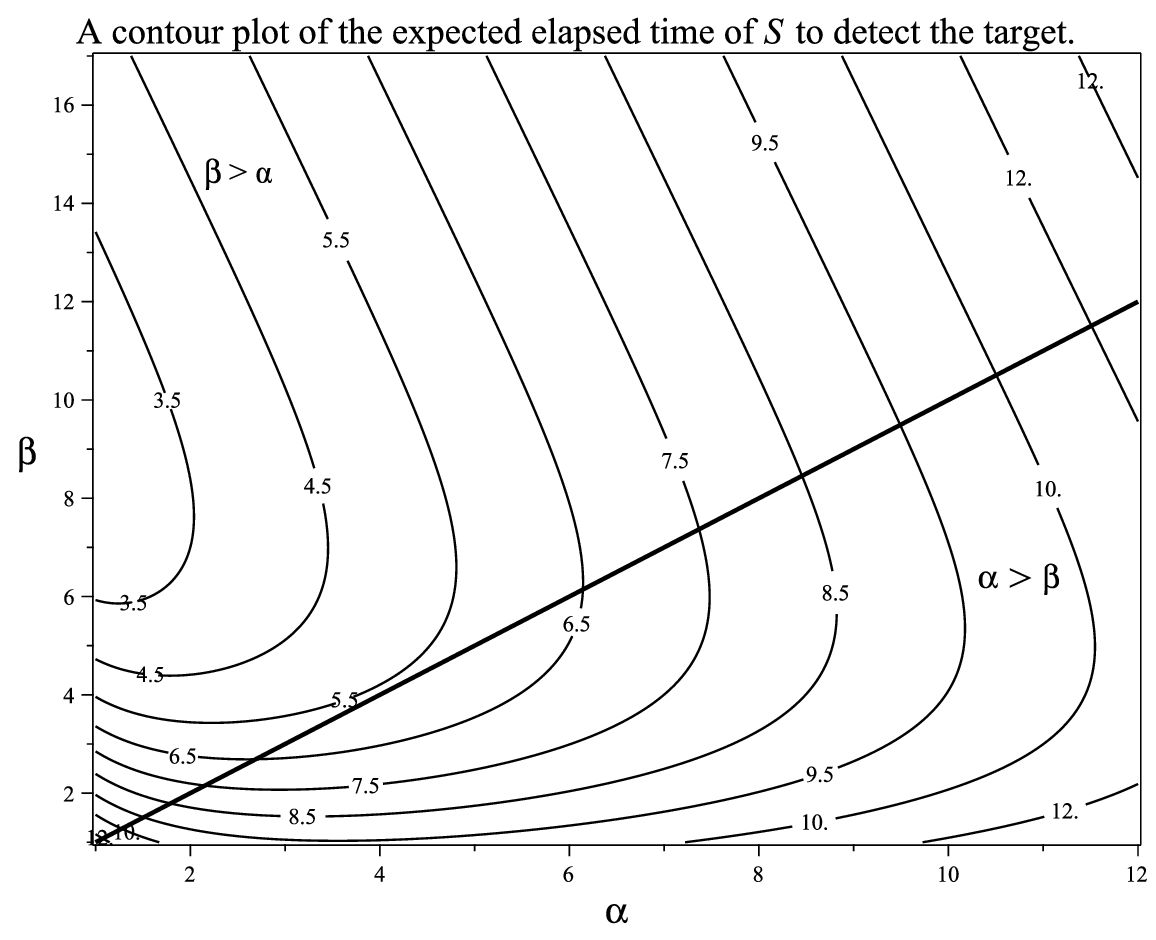}\\
  \vspace{0.0cm}\caption{\footnotesize A contour plot of the expected time of a target follows gamma distribution as a function of the parameters $\alpha$ and $\beta$.}\label{Fig11}
\end{figure}

A search problem consists of two parts; the first is the behavior of the hidden target, i.e., the distribution that the target follows. The second part is the mechanism of the searcher. This paper proposes applying several truncations within a bounded interval based on exploiting the available information about the target to enhance search efficiency. In this section, we present five examples to discuss the properties of the proposed technique. The first three examples are dedicated to discussing the nature of the truncation and the resultant truncated distributions. The first example demonstrates the effect of commensurate truncation on various distributions and their CDFs. The impact of uneven truncation on some distributions is discussed in Example 2. Next, the symmetric truncation on symmetric distributions is presented in Example 3. The last two examples are dedicated to discussing the expected elapsed time for each case and optimizing these expectations. The simulation is carried out using the Maple 2019 program on a 2.7GHz Xeon processor. \\
\begin{exm}\label{exmp1}
Here we study the effects of the commensurate truncation on symmetric distributions such as normal and Cauchy distributions; also apply it on asymmetric distributions, e.g., skew-normal distributions, and a truncation with $N$-deleted subintervals is applied on gamma distribution, their PDFs and CDFs are given Tables \ref{Tab1}, and \ref{Tab2}. Figures \ref{Fig4}(a)-(c) represent the original normal, Cauchy, and skew-normal distributions and their corresponding truncated distributions under the commensurate case when $N=2$, and endpoints a=-20, b=30. The values of their parameters are given in Table \ref{Tab3} and the used truncation mesh points are listed in Table \ref{Tabprmexm1}. In light of equation \eqref{Gfuns1}, the commensurate truncated CDFs of these distributions are plotted in Figures \ref{Fig5}(a)-(c). It has been observed that they tend to $0$, $1$ as $x$ approaches the endpoints $a$ and $b$, i.e., satisfy the intrinsic property of the cumulative distribution function. The PDF of the truncated gamma distribution and the associated CDF that is calculated using equation  \eqref{Ufuns3} are displayed in Figures \ref{Fig6}(a), and \ref{Fig6}(b). The area of shaded regions in Figures \ref{Fig4}(a)-(c), and \ref{Fig6}(a), are 1, which represents the total probability in the considered domain. Moreover, the resultant truncated PDFs of the symmetric distributions are not symmetric.
\begin{center}
\begin{table}
\begin{center}
  \begin{tabular}{p{5.5cm}p{6.0cm}}
    \hline
    \textbf{Distribution name} &  \textbf{The values of its parameters}   \\
    \hline
    Normal distribution & $\sigma=4.53,~~\mu=0$, \\
    Cauchy distribution & $\tilde{x}=1.38,~~c=2.76$,\\
    Skew normal distribution & $\varrho=-3.37$, $\eta=-2.3$, $\varpi=1.6$, \\
    Gamma distribution & $\kappa=3.15$, ~$\theta=1.27$.\\
    \hline
  \end{tabular}
\end{center}
 \caption{The values of distributions parameters used in Example \ref{exmp1}.}\label{Tab3}
\end{table}
\end{center}
\begin{center}
\begin{table}
\begin{center}
  \begin{tabular}{p{1.0cm}p{2.0cm}p{2.0cm}p{2.0cm}p{2.0cm}}
    \hline
    $i$ &  $\vartheta_{i}$ &  $\zeta_{i}$ & $\alpha_{i}$ & $\beta_{i}$  \\
    \hline
    $1$& $-6$ & $-4$  & $2$  & $7$  \\
    $2$& $-15$ & $-10$  & $11$  & $17$  \\
    \hline
  \end{tabular}
\end{center}
 \caption{The values of truncation mesh points used in commensurate case.}\label{Tabprmexm1}
\end{table}
\end{center}
\end{exm}
\begin{exm}\label{exmp2}
A truncation is said to uneven when the number of deleted subintervals of the left and right parts are not equal. Here we apply an uneven truncation on normal, Cauchy and skew-normal distributions. Figures \ref{Fig7}(a)-(c) depict The resultant truncated PDFs when $M=3$, and $N=2$, with distribution parameters are given in Table \ref{Tab4}, the truncation mesh points are given in Table \ref{Tabprmexm2}, and the endpoints are $a=-30$, and $b=20$. Based on equation \eqref{Hfuns2}, the corresponding CDFs are plotted in Figures \ref{Fig8}(a)-(c).
\begin{center}
\begin{table}
\begin{center}
  \begin{tabular}{p{5.5cm}p{6.0cm}}
    \hline
    \textbf{Distribution name} &  \textbf{The values of its parameters}   \\
    \hline
    Normal distribution & $\sigma=3.6$,~~$\mu=2.5$, \\
    Cauchy distribution & $\tilde{x}=4.6$,~~$c=2.5$,\\
    Skew normal distribution & $\varrho=5.48$, $\eta=-2.13$, $\varpi=2.6$. \\
    \hline
  \end{tabular}
\end{center}
 \caption{The values of distributions parameters of Example \ref{exmp2}.}\label{Tab4}
\end{table}
\end{center}
\begin{center}
\begin{table}
\begin{center}
  \begin{tabular}{p{1.0cm}p{2.0cm}p{2.0cm}p{2.0cm}p{2.0cm}}
    \hline
    $i$ &  $\vartheta_{i}$ &  $\zeta_{i}$ & $\alpha_{i}$ & $\beta_{i}$  \\
    \hline
    $1$& $-6$ & $-4$  & $3$  & $6$  \\
    $2$& $-17$ & $-13$  & $13$  & $17$  \\
    $3$& $-23$ & $-20$  & --  &--  \\
    \hline
  \end{tabular}
\end{center}
 \caption{The values of mesh points used in uneven truncation.}\label{Tabprmexm2}
\end{table}
\end{center}
\end{exm}
\begin{exm}\label{exmp3}
So far, the resultant distributions after applying truncations are asymmetric, as shown in examples \ref{exmp1} and  \ref{exmp2}. In this example, we discuss the circumstances that result in symmetric truncated distribution. The symmetric truncation requires that the number of deleted subintervals of the left and right parts is equal ( $M=N$). Also, the truncation mesh points in the left part of the origin are equal to the negative of the corresponding mesh points of the right part, i.e., $\zeta_{i}=-\alpha_{i}$, $\vartheta_{i}=-\beta_{i}$ for all values of $i\in\{1,2,\ldots,N\}$. Finally, the value of the endpoints must satisfy $a=-b$. Here we have $N=2$; based on these conditions, the following truncation mesh points are selected; $(\alpha_{1},\beta_{1})=(4,8)$, and $(\alpha_{2},\beta_{2})=(13,17)$, the endpoints are $a=-20$, and $b=20$. Next, the original distribution must be symmetric, so we apply symmetric truncation on a normal distribution with parameters $\sigma=4.53$, and $\mu=0$, and on Cauchy distribution with parameters $\tilde{x}=0$, and $c=2.76$. Figure \ref{Fig9} is drawn in light of these conditions and values.
\end{exm}
\begin{exm}\label{exmp4}
In this example, we discuss the expected time values that the searchers $S_{1}$ and $S_{2}$ require to find the lost target. The searcher $S_{1}$ sweeps the right part of the line, while $S_{2}$ scans the left part. They move with velocity $v_{1}$ in the search zones, whereas their velocities in the skipping zones are $v_{2}$. Here we choose the value of $v_{2}=5v_{1}$, and $v_{1}$ is a unit velocity, i.e., $v_{1}=1$. The expected elapsed time is investigated from two aspects; the first, we calculate $\mathbb{E}(\tau_{m})$, given the values of the truncation mesh points. Second, we consider the value of the truncation mesh points is unknown. We analyze and optimize them based on the given distribution to obtain an optimum expectation value. For the first case, the expected time values $\left\{\mathbb{E}(\tau_{m})\right\}_{m=0}^{\hat{N}+1}$ are calculated considering the commensurate and uneven truncations that applied on distributions of examples 1, and 2. Based on equations \eqref{expcttmeq1} and \eqref{expcttmeq2}, the sequences of expected elapsed times are calculated and listed in Table \ref{Tab7} of the normal, Cauchy and skew-normal distributions for commensurate truncation with $N=2$, and uneven truncation with $M=3$, and $N=2$, respectively. Also, we calculate the expected elapsed time of sweeping the left and right parts of the search zone without applying any truncation within it; their values are given in Table \ref{Tab8}. Note that $\mathbb{E}(\tau_{0})$, and  represents the expected elapsed time when the hidden object lies in the interval $[a,\vartheta_{\tilde{M}}]$, and $\mathbb{E}(\tau_{\hat{N}+1})$ is the time expected value if the lost target lies in the interval $[\beta_{N},b]$. It has been observed that applying truncations within the search zone reduces the search cost to find the lost target.\\
\begin{center}
\begin{table}
\begin{center}
\begin{tabular}{cc||c||c||c|}
\cline{2-5}
&\multicolumn{4}{| c| }{\textbf{Commensurate truncation} $(N=2)$}\\
\cline{2-5}
&\multicolumn{1}{|c|}{$\scriptsize{\mathbb{E}(\tau_{m})}$} & \textbf{Normal distribution } & \textbf{Cauchy distribution}& \textbf{Skew-normal distribution}\\
\cline{2-5}
&\multicolumn{1}{|c|}{ $\mathbb{E}(\tau_{0})$ }& $3.94262$ &$3.25040$&$5.40$\\
& \multicolumn{1}{|c|}{$\mathbb{E}(\tau_{1})$} & $2.93889$ &$2.13912$&$4.40$\\
& \multicolumn{1}{|c|}{$\mathbb{E}(\tau_{2})$ }& $2.02496$ &$1.462140$&$3.88673$\\
&\multicolumn{1}{|c|}{ $\mathbb{E}(\tau_{3})$ }& $0.55463$ &$0.790499$&$2.6e$-8\\
 \parbox[t]{2mm}{\multirow{1}{*}{\rotatebox[origin=c]{90}{Example \ref{exmp1}}}}
& \multicolumn{1}{|c|}{$\mathbb{E}(\tau_{4})$} & $1.90292$ &$2.19950$&1.0\\
& \multicolumn{1}{|c|}{$\mathbb{E}(\tau_{5})$} & $3.1048$ &$3.99056$&2.2\\
\cline{2-5}
&\multicolumn{4}{| c| }{\textbf{Uneven truncation} $(M=3,~~N=2)$}\\
\cline{2-5}
&\multicolumn{1}{|c|}{ $\mathbb{E}(\tau_{0})$ }& 3.99642 &3.17277&6.11292\\
& \multicolumn{1}{|c|}{$\mathbb{E}(\tau_{1})$} & 3.39642&2.49555&5.51292\\
& \multicolumn{1}{|c|}{$\mathbb{E}(\tau_{2})$ }& 1.69452 &1.30486&2.86452 \\
&\multicolumn{1}{|c|}{ $\mathbb{E}(\tau_{3})$ }& 1.20252 &0.520927&2.46452\\
 \parbox[t]{2mm}{\multirow{1}{*}{\rotatebox[origin=c]{90}{Example \ref{exmp2}}}}
& \multicolumn{1}{|c|}{$\mathbb{E}(\tau_{4})$} & 1.34945 &0.915580&1.14605\\
& \multicolumn{1}{|c|}{$\mathbb{E}(\tau_{5})$} & 3.60403 &4.78605&1.75902\\
& \multicolumn{1}{|c|}{$\mathbb{E}(\tau_{6})$} & 4.40415 &5.65515&2.55902\\
\cline{2-5}
\end{tabular}
\end{center}
\caption{\footnotesize{The expected elapsed time to find a hidden target follows various distributions of examples \ref{exmp1} and \ref{exmp2} considering the commensurate and uneven truncations.}}\label{Tab7}
\end{table}
\end{center}
\begin{center}
\begin{table}
\begin{center}
\begin{tabular}{c|c|c|c|}
  \hline
  \multicolumn{4}{|c|}{\textbf{Example} \ref{exmp1}}\\
  \hline
\multicolumn{1}{|c|}{\textbf{  The expectation}} & \small\textbf{Normal distribution } & \small\textbf{Cauchy distribution}  &\small\textbf{Skew-normal distribution}  \\
  \hline
   \multicolumn{1}{|c|}{of the left part} & 9.99995 &6.71061 &20.0\\
 \multicolumn{1}{|c|}{ of the right part} & 15.0    &19.9341&2.8$e$-7\\
  \hline
    \multicolumn{4}{|c|}{\textbf{Example} \ref{exmp2}}\\
  \hline
   \multicolumn{1}{|c|}{of the left part} & 7.31106&4.390810&17.62036\\
 \multicolumn{1}{|c|}{ of the right part} & 15.1260 &17.07280&8.25310\\
  \hline
\end{tabular}
\end{center}
\caption{\footnotesize{The expected elapsed time to find a hidden target follows various distributions of examples \ref{exmp1} and \ref{exmp2} without applying any truncation within the search zone.}}\label{Tab8}
\end{table}
\end{center}
Next, we illustrate the idea of determining the value of truncation mesh points that improve the value of the expected elapsed time. Here we consider two cases; a commensurate truncation with $N=1$ deleted interval in each part that applied on a normal distribution. The second case; applying a truncation with one deleted subinterval on gamma distribution with parameters given in Table \ref{Tab3} of Example \ref{exmp1}. By substituting $N=1$ into equation \eqref{expcttmeq1}, we obtain $\mathbb{E}(\tau_{0})$ as a function of $\vartheta$, and $\zeta$, while $\mathbb{E}(\tau_{3})$ is a function of $\alpha$, and $\beta$ for the normal distribution. In the case of the gamma distribution, we put $N=1$ into equation \eqref{expcttmeq3}, then we have $\mathbb{E}(\tau_{1})$ represents the total time to sweep the search zone, which it is a function in the variable $\alpha$ and $\beta$. The indices are omitted since we have one of each variable. After that, we plot contour figures of these expectations as revealed in Figure \ref{Fig10} for the normal distribution and  Figure \ref{Fig11} for a target following gamma distribution. The straight lines in Figures \ref{Fig10} and \ref{Fig11} represent the value of the expected elapsed time when $\vartheta=\zeta$, and $\alpha=\beta$. Based on the proposed model, we have $\vartheta<\zeta$ and $\alpha<\beta$ consequently, we can choose suitable values for these variables from the region above the line that minimizes the expected values.\\
\end{exm}
\begin{exm}\label{exmp5}
Finally, we discuss optimizing the expected elapsed time values of the searchers $S_{1}$ and $S_{2}$ under the commensurate and uneven cases. Based on pondering the whole problem, we have two optimization cases; the first case is to optimize the required time by one of the two searchers, i.e., applying the optimization on either $\mathbb{E}(\tau_{0})$ for searcher $S_{2}$, or $\mathbb{E}(\tau_{\breve{N}})$ for the searcher $S_{1}$. The second case is to optimize the elapsed time to sweep the whole truncated domain considering the truncations on the negative and positive parts of the real line together. Here, we consider the second case to optimize the total expected elapsed time of the search process. In other words, we optimize the function $\hat{\mathbb{E}}(\tau_{0},~\tau_{\breve{N}})$, such that
$$\hat{\mathbb{E}}(\tau_{0},~\tau_{\breve{N}})=\hat{\mathbb{E}}(\tau_{0})+\hat{\mathbb{E}}(\tau_{\breve{N}}),$$
which obtained by adding equations \eqref{mdfEsym1}and \eqref{mdfEsym2}. We apply the optimization algorithm \ref{Algthm1} on $\hat{\mathbb{E}}(\tau_{0},~\tau_{\breve{N}})$, for the commensurate case when $N=2$, $3$, and uneven case when $M=3$, $N=2$. Table \ref{TabSpOpt1} reports the optimized values of $\{(\vartheta_{i},\zeta_{i})\}_{i=1}^{\tilde{M}}$, and $\{(\alpha_{i},\beta_{i})\}_{i=1}^{N}$ that minimize the total expected time value to sweep the whole search domain to find a lost object considering various distributions such as normal, Cauchy, and skew-normal distributions. The values of the distribution parameters under the commensurate truncation are listed in Table \ref{Tab3}, and their values under the uneven truncation are given in Table \ref{Tab4}. Also, the optimized time expectations and the corresponding CPU times are listed in Table \ref{TabSpOpt1}.
\begin{center}
\tiny\begin{table}
\begin{center}
\begin{adjustbox}{angle=90}
\tiny\begin{tabular}{c|c|c|c|c|c||c|c|c|c||c|c|c|c|}
  \cline{3-14}
\multicolumn{2}{ c }{}  &\multicolumn{12}{| c|}{\textbf{Commensurate truncation}}\\
\cline{3-14}
\multicolumn{2}{ c }{} & \multicolumn{4}{| c||}{\textbf{Normal distribution}}&\multicolumn{4}{ c|| }{\textbf{Cauchy distribution}}& \multicolumn{4}{ c| }{\textbf{Skew-normal distribution}} \\
\cline{2-14}
&\multicolumn{1}{ c| }{$i$} & $\vartheta_{i}^{\ast}$ & $\zeta_{i}^{\ast}$ & $\alpha_{i}^{\ast}$ & $\beta_{i}^{\ast}$ &$\vartheta_{i}^{\ast}$ & $\zeta_{i}^{\ast}$ & $\alpha_{i}^{\ast}$ & $\beta_{i}^{\ast}$& $\vartheta_{i}^{\ast}$ & $\zeta_{i}^{\ast}$ & $\alpha_{i}^{\ast}$ & $\beta_{i}^{\ast}$ \\
\cline{2-14}
&\multicolumn{1}{ c| }{$1$} & -5.2816 & -4.5546 & 3.6613 & 5.3652 & -5.2681 & -4.5686 & 3.7372 & 5.3266 & -5.1898 & -4.2622 & 3.7333& 5.2667   \\
 \parbox[t]{2mm}{\multirow{1}{*}{\rotatebox[origin=c]{90}{$N=2$}}}
&\multicolumn{1}{|c| }{$2$} & -13.2232 & -11.6417 & 12.9763 & 14.8908 & -13.2643 & -11.6822 & 13.0057 & 14.9978 & -13.2389 & -11.6778 & 13.0667 & 14.9333   \\
\cline{2-14}
& \multicolumn{3}{ c}{$\mathbb{E}^{\ast}(\tau_{0},~\tau_{\breve{N}})$} & \multicolumn{2}{| c||}{$6.1034$} & \multicolumn{4}{ c||}{$6.4001$}& \multicolumn{4}{ c|}{5.6254}   \\
\cline{2-14}
& \multicolumn{3}{ c}{\tiny\textbf{CPU time (sec)}} & \multicolumn{2}{ c||}{$0.860$} & \multicolumn{4}{ c||}{$0.703$} & \multicolumn{4}{ c|}{$1.407$}   \\
\cline{2-14}
&\multicolumn{1}{ c| }{$i$} & $\vartheta_{i}^{\ast}$ & $\zeta_{i}^{\ast}$ & $\alpha_{i}^{\ast}$ & $\beta_{i}^{\ast}$ &$\vartheta_{i}^{\ast}$ & $\zeta_{i}^{\ast}$ & $\alpha_{i}^{\ast}$ & $\beta_{i}^{\ast}$& $\vartheta_{i}^{\ast}$ & $\zeta_{i}^{\ast}$ & $\alpha_{i}^{\ast}$ & $\beta_{i}^{\ast}$ \\
\cline{2-14}
&\multicolumn{1}{ c| }{$1$} & $-8.0892$  & $-7.3834$  & $4.2824$  & $5.6125$  & $-8.1654$ & $-7.5199$   & $4.2056$ & $5.5873$ & $-7.9318$ & $-7.0661$ & $4.4000$ & $5.6000$   \\
&\multicolumn{1}{ c| }{$2$} & $-12.9223$ & $-12.0460$ & $12.3246$ & $13.5645$ & $-12.9304$ & $-12.0553$ & $12.3289$ & $13.5919$ & $-12.9332$ & $-12.0664$ & $12.4000$ & $13.6000$   \\
 \parbox[t]{2mm}{\multirow{1}{*}{\rotatebox[origin=c]{90}{$N=3$}}}
&\multicolumn{1}{| c| }{$3$} & $-17.2659$ & $-16.7319$ & $22.3990$ & $23.5995$ & $-17.2626$ & $-16.7268$ & $22.3762$ & $23.5940$ & $-17.2667$ & $-16.7333$ & $22.4000$ & $23.6000$   \\
\cline{2-14}
& \multicolumn{3}{ c|}{$\mathbb{E}^{\ast}(\tau_{0},~\tau_{\breve{N}})$} & \multicolumn{2}{c||}{$7.2462$} & \multicolumn{4}{ c||}{$6.8690$}& \multicolumn{4}{ c|}{$8.2380$}   \\
\cline{2-14}
& \multicolumn{3}{ c|}{\tiny\textbf{CPU time (sec)}} & \multicolumn{2}{ c||}{$1.672$} & \multicolumn{4}{ c||}{$1.265$} & \multicolumn{4}{ c|}{$5.155$}   \\
\cline{2-14}
&\multicolumn{1}{ c }{}  &\multicolumn{12}{ c|}{\textbf{uneven truncation} $M=3$, $N=2$}\\
\cline{2-14}
&\multicolumn{1}{ c| }{$i$} & $\vartheta_{i}^{\ast}$ & $\zeta_{i}^{\ast}$ & $\alpha_{i}^{\ast}$ & $\beta_{i}^{\ast}$ &$\vartheta_{i}^{\ast}$ & $\zeta_{i}^{\ast}$ & $\alpha_{i}^{\ast}$ & $\beta_{i}^{\ast}$& $\vartheta_{i}^{\ast}$ & $\zeta_{i}^{\ast}$ & $\alpha_{i}^{\ast}$ & $\beta_{i}^{\ast}$ \\
\cline{2-14}
&\multicolumn{1}{ c| }{$1$} & $-5.2316$  & $-4.6001$  & $4.0543$  & $5.2179$  & $-5.2903$   & $-4.7060$  & $4.2885$ & $5.4200$   & $-5.3177$ & $-4.3357$ & $4.0667$ & $4.9333$   \\
&\multicolumn{1}{ c| }{$2$} & $-15.5940$ & $-14.3880$ & $14.2291$ & $15.5145$ & $-15.5870$  & $-14.3656$ & $14.1120$ & $15.4616$ & $-15.5769$ & $-14.3538$ & $14.4000$ & $15.6000$  \\
&\multicolumn{1}{ c| }{$3$} & $-21.9333$ & $-21.0667$ & --        & --        & $-21..9378$ & $-21.0642$ & --        & --        & $-21.9333$ & $-21.0667$ & -- & --   \\
  \cline{2-14}
& \multicolumn{3}{ c|}{$\mathbb{E}^{\ast}(\tau_{0},~\tau_{\breve{N}})$} & \multicolumn{2}{ c||}{$6.5724$} & \multicolumn{4}{ c||}{$7.4671$}& \multicolumn{4}{ c|}{$5.6890$}   \\
\cline{2-14}
& \multicolumn{3}{ c|}{\tiny\textbf{CPU time (sec)}} & \multicolumn{2}{ c||}{$1.500$} & \multicolumn{4}{ c||}{$1.063$} & \multicolumn{4}{ c|}{$3.235$}   \\
\cline{2-14}
\end{tabular}
\end{adjustbox}
\end{center}
 \caption{\tiny The values of the optimized mesh truncation points of Example \ref{exmp5} for several distributions under commensurate and uneven cases, the optimized expected time values, and the CPU time.}\label{TabSpOpt1}
\end{table}
\end{center}
\end{exm}
\newpage
\large\section{Conclusion and future work}
Sweeping a search zone has several deleted regions with two searchers in one dimension leads to minimizing the expected elapsed time to find the lost target. Deleting multiple regions results in three categories of truncations based on their numbers in the left and right parts of the observer center and on the values of the mesh truncation points. We define the corresponding truncated PDF and its CDF by considering the truncation classification. We considered the lost target follows some distributions, which can be symmetric such as normal and Cauchy distributions, or asymmetric such as skew-normal distribution. Most published methods are not concerned with truncating the intervals from the target distribution where the target is poor probably to be found (they are truncated based on the available information about the target's position), so the expected cost increases. Next, the expected elapsed time is introduced for each case and optimized based on Newton's method to reduce the required time to find the lost target. Five examples are executed to illustrate the truncated distributions' behavior and calculate the time needed to find the missing target.\\\\
In future research, it is possible to apply our model to detect a target located on one of the disjointed real lines. In addition, the multiple truncation definition for a multivariate distribution can be presented to test our model's ability to detect lost targets.
\section*{Data Availability}
The data used to support the findings of this study are included within the article.
\section*{Acknowledgement}
This research has been funded by Princess Nourah bint Abdulrahman University and Researchers Supporting Project number (PNURSP2023R346), Princess Nourah bint Abdulrahman University, Riyadh, Saudi Arabia.
\section*{Conflict of interest}
The authors declare that they have no conflict of interest.
\newpage
\bibliography{TruncationRefernces.bib}
\bibliographystyle{ieeetr}

\end{document}